\documentclass[11pt,a4paper]{article}
\usepackage[margin=1in]{geometry}
\usepackage[T1]{fontenc}
\usepackage{lmodern}
\usepackage{amsmath,amssymb,amsthm,mathtools}
\usepackage{microtype}
\usepackage{enumitem}
\usepackage{hyperref}
\usepackage{xcolor}
\usepackage{comment}
\usepackage{marginnote} 
\usepackage{natbib}
\usepackage{url}
\newcommand{\revnote}[1]{\leavevmode\marginnote{\raggedright\fontsize{5.6}{6.3}\selectfont\color{red}#1}}
\newcommand{\redrev}[2]{\revnote{#1}{\color{red}#2}}
\hypersetup{colorlinks=true,linkcolor=blue,citecolor=blue,urlcolor=blue}

\newtheorem{theorem}{Theorem}[section]
\newtheorem{proposition}[theorem]{Proposition}
\newtheorem{lemma}[theorem]{Lemma}
\newtheorem{corollary}[theorem]{Corollary}
\theoremstyle{remark}
\newtheorem{remark}[theorem]{Remark}

\newcommand{\R}{\mathbb R}

\newcommand{\norm}[2][]{\left\|#2\right\|_{#1}}
\newcommand{\ip}[2]{\left\langle #1,#2\right\rangle}
\newcommand{\dd}{\mathrm d}
\newcommand{\D}{\partial_x}
\newcommand{\M}{\mathrm M}
\newcommand{\F}{\mathcal F}

\newcommand\warning[1]{#1} 

\makeatletter
\newcommand{\multilabel}[1]{\ltx@label{#1}}
\makeatother

\title{Quantitative constructive representation of the KdV flow in $H^3$ 
via $L^2$-energy-preserving centered-difference Chernoff approximations and Picard iteration}

\author{
Comlan Edmond Koudjinan\\[4pt]
\small  University of Toronto, Toronto, Ontario, Canada\\[4pt]
\small \href{mailto:Koudjinanedmond@gmail.com}
{\texttt{Koudjinanedmond@gmail.com}}
}
\date{}

\date{\today}

\begin{document}
\maketitle

\begin{abstract}
We give a constructive representation of the solution of the KdV equation $u_t=-u_{xxx}+6uu_x$ on $\mathbb R$ for real initial data in $H^3(\mathbb R)$. The approximations couple an exact, unitary centered-difference Chernoff product for the Airy group with a skew-symmetric discretization of the quadratic term, and are computed by Picard iteration; they preserve the $L^2$ norm exactly, as does the limit. A two-mesh Cauchy argument 
shows that the whole sequence converges in $C([0,T];H^2)$, at rate $O(n^{-\frac{2(s-2)}5})$ for $H^s$ data ($3\leq s\leq 7$), and in $L^2$ at rate $O(n^{-6/5})$ for $H^3$ data. 
This scheme, to the best of our knowledge, is novel and is of independent interest. As a by-product we obtain existence and uniqueness of strong solutions, Lipschitz dependence on the data in the $H^2$ norm, uniformly on bounded subsets of $H^3$ and compact time intervals, and joint locally Lipschitz dependence of $(t,f)\mapsto S(t)f$ into $L^2$; the classical conservation laws extend the flow globally in time.
\end{abstract}

\noindent\textbf{MSC 2020.}
35Q53; 35B30, 35A35, 47H20, 65M06, 65M12, 65M15.

\medskip
\noindent\textbf{Keywords.}
KdV equation; Chernoff product formula;
centered differences; $L^2$-norm-preserving approximation;
Picard iteration; convergence rates;
constructive representation.
\section*{Introduction}

We consider the Cauchy problem for the Korteweg--de Vries (KdV)
equation on the real line,
\begin{equation}\label{eq:KdV.init}
u_t+u_{xxx}-6uu_x=0,\qquad u(0)=u_0\in H^3(\R;\R),
\end{equation}
derived by \citet{KortewegDeVries1895} as a model for long, weakly
nonlinear shallow-water waves. Its solution is obtained as the limit
of explicit $L^2$-norm-preserving approximations, built from centered
translations, a Chernoff product for the Airy group and Picard
iteration, and we quantify the rate of convergence in Sobolev norms.

\paragraph{Fixed-point methods.}
For regular data, \eqref{eq:KdV.init} is well posed by energy methods
combined with approximation or regularization
\citep{Sjoberg1970,BonaSmith1975,Kato1983}. Most of the low-regularity
theory instead applies the Banach fixed-point theorem to the Duhamel
formulation
\begin{equation}\label{eq:intro-Duhamel}
u(t)=e^{tL}u_0+\int_0^t e^{(t-s)L}\,6u(s)\D u(s)\,\dd s,
\qquad L=-\D^3,
\end{equation}
where $(e^{tL})_{t\in\R}$ is the Airy group. The loss of derivative in the nonlinearity
 prevents a naive application the contraction principle in
$C([0,T];H^s)$; one way of recovering this loss is to use the dispersion property, either
through local smoothing and maximal function estimates
\citep{KenigPonceVega1991,KenigPonceVega1993} or in the Fourier
restriction spaces of \citet{Bourgain1993KdV},
\begin{equation}\label{eq:intro-Xsb}
\|u\|_{X^{s,b}}
=\bigl\|\langle\xi\rangle^{s}\langle\tau-\xi^{3}\rangle^{b}\,
\widetilde u(\tau,\xi)\bigr\|_{L^2_{\tau,\xi}},
\end{equation}
where $\widetilde u$ is the space--time Fourier transform. 
The index $s$
measures spatial regularity, while $b$ measures the time regularity of
the profile $e^{-tL}u(t)$, equivalently the concentration of
$\widetilde u$ near the surface $\tau=\xi^3$.
The
bilinear estimate of \citet{KenigPonceVega1996} in these spaces yields
local well-posedness in $H^s(\R)$ for $s>-3/4$; see
\citep{CKSTT01,Guo2009,Kishimoto2009} for global solutions and the
endpoint, and \citet{Tao2006} for an exposition. The sharp result in
$H^{-1}(\R)$ \citep{KillipVisan2019} relies instead on complete
integrability.

The contraction principle yields the solution as the limit of Picard
iterates. 
Each
iterate, however, involves the Airy group and the Duhamel
integral, the contraction is proved in a space--time Fourier norm:
it is an existence mechanism rather than a
pratical way to find the solution.

\paragraph{Constructive approximation and rates.}
Solutions have also been obtained as limits of numerical schemes:
semi-discrete finite differences for the periodic problem
\citep{Sjoberg1970}, and fully discrete finite-difference and
Galerkin-type schemes converging to strong solutions for $H^3$ data
and to weak solutions for $L^2$ data
\citep{HoldenKoleyRisebro2015,DuttaKoleyRisebro2015}. 
However, these constructions do not quantify the convergence; in
\citet{HoldenKoleyRisebro2015} and \citet{DuttaKoleyRisebro2015}, for
instance, it follows from uniform bounds and compactness. 
We refer also the reader, for operator splitting
\citep{HKRT2011,HoldenLubichRisebro2013}, finite-difference schemes
\citep{CourtesLagoutiereRousset2020} and exponential-type integrators
\citep{HofmanovaSchratz2017}. 
As far as we know, we are not aware of any construction that
 produces the solution of \eqref{eq:KdV.init} as a limit of
approximations and establishes an explicit rate for this convergence. One of the goals of this work is to fill this gap.

\paragraph{Chernoff approximations.}
Chernoff's theorem \citep{Chernoff1968} states that
$\lim_{n\to \infty}F(t/n)^nf\to e^{tA}f$, uniformly for $t$ in compact intervals,
whenever $F:[0,\infty)\to\mathcal L(X)$ is strongly continuous, $F(0)=I$, and $F'(0)=A$ on a core of the generator
$A$, under the stability property $\|e^{tA}\|\le Me^{\omega t}$ and $\|F(t)^k\|\le Me^{\omega kt}$, $t\geq 0$ and $k\in\mathbb N$ \citep[see also][]{EngelNagel2000,Butko2020}. Remizov and
coauthors developed this into a method for explicit solution formulas
for linear evolution equations with variable coefficients, using
Chernoff functions built from translation operators
\citep{Remizov2016MN,Remizov2018AMC,Remizov2019JMP} and quasi-Feynman
formulas for Schr\"odinger groups \citep{Remizov2016JFA}; rates of
convergence, and Chernoff functions with faster convergence, are
derived in \citet{GalkinRemizov2025} and \citet{VedeninRemizov2026}.
For linear evolutions, centered differences and Richardson
extrapolation give arbitrarily high rate of convergence and, under
analytic regularity, super-fast rates
\citep{Koudjinan2026SuperfastChernoff}.

\paragraph{The present construction.}
%
In this paper, we approximate the linear part, the Airy group generated by the
skew-adjoint operator $L=-\D^3$. For, we use the unitary Chernoff function
$\mathcal S_2(h)=\exp(hQ_h)$, where $Q_h=-\delta_h^3+\frac{h^2}{2}\delta_h^5$
is built from the centered difference
$\delta_h=(T_h-T_{-h})/(2h)$, so that the linear propagator is the
Chernoff product $U_n(t,0)=e^{B_n(t)}=\mathcal S_2(t/n)^n$
(Theorem~\ref{thm:airy-product}, where the convergence is verified
directly on the Fourier side). The quadratic term is discretized by the
split form $\F_h(f)=2f\delta_hf+2\delta_h(f^2)$ of \eqref{eq:Fh}, which
is skew-symmetric in the sense that $\F_h(q)=\mathcal K_h(q)q$ with
$\mathcal K_h(q)^*=-\mathcal K_h(q)$ \eqref{eq:K}; its three-point
average is of the type used in the early finite-difference treatment
of KdV by \citet{ZabuskyKruskal1965}. The
approximations $p_{n,\varepsilon}$ are computed by the Picard iteration
\eqref{eq:Picard-constructive}. The resulting evolutions are
nonautonomous equations on $\R$ in continuous time.

The contraction principle enters here differently from the fixed-point
methods above. For a fixed positive nonlinear mesh,
$\|\delta_h\|=|h|^{-1}$ on every $H^r(\R)$, so the Duhamel map
\eqref{eq:Picard-constructive} is a contraction on $C(I_j;H^3)$ for
intervals of length $O(\varepsilon/(n\|f\|_{H^3}))$
(Lemma~\ref{lem:positive-mesh}); no dispersive function space is
needed, and each iterate is computed from centered translations and
explicit Fourier multipliers. These intervals shrink as $n\to\infty$,
so the contraction provides no control uniform in $n$. Mesh-uniform
control comes instead from the skew structure: $L^2$-energy conservation 
and an $H^3$ energy estimate for $\F_h$ with constants independent of
$h$ (Lemma~\ref{lem:H3}) yield
$\|p_{n,\varepsilon}(t)\|_{H^3}\le4\|f\|_{H^3}$ on the common interval
$[0,T_*]$, $T_*=c_*/\sqrt{\|f\|_2\|f\|_{H^3}}$
(Proposition~\ref{prop:H3apriori} and Lemma~\ref{lem:positive-mesh}).
Convergence then follows from a two-mesh argument. The difference
$q_n-q_m$ solves an equation whose nonlinear part is controlled by a
one-sided $H^2$ estimate without derivative loss
(Lemma~\ref{lem:H2onesided}), while the mismatch of the generators and
nonlinearities at the meshes $t/n$ and $t/m$ is of order
$t^2|n^{-2}-m^{-2}|$. For $H^7$ data this gives
$\|q_n-q_m\|_{C([0,T];H^2)}\le C|n^{-2}-m^{-2}|$; for $H^3$ data a
density argument shows that the whole sequence is Cauchy in
$C([0,T];H^2)$ 
(Theorem~\ref{thm:order2}).

\paragraph{Main results.}
For real $0\ne f\in H^3(\R)$, $T=T_*(f)$ and $0<\varepsilon\le T$,
Theorem~\ref{prop:explicit-formula} gives
$\|p_{n,\varepsilon}(t)\|_2=\|f\|_2$,
$p_{n,\varepsilon}\to S(\cdot)f$ as $n\to \infty$ in $C([0,T];H^2)$, and
\[
\|p_{n,\varepsilon}-S(\cdot)f\|_{C([0,T];L^2)}\le C_{T,f}\,n^{-6/5},
\qquad
\|p_{n,\varepsilon}-S(\cdot)f\|_{C([0,T];H^2)}
\le C_{T,f,s}\,n^{-\frac{2(s-2)}{5}}
\]
for $f\in H^s(\R)$, $3\le s\le7$; the convergence is quadratic for
$H^7$ data. As a by-product,
$S(\cdot)f\in C([0,T];H^3)\cap C^1([0,T];L^2)$ is the unique solution
in $C([0,T];H^2)\cap C^1([0,T];H^{-1})$
(Lemma~\ref{lem:kdv-H2-uniqueness}). By the classical conservation
laws the flow is global, Lipschitz in the data with respect to the
$H^2$ norm on bounded subsets of $H^3$, H\"older continuous of
exponent $1/3$ in time with values in $H^2$, and jointly locally
Lipschitz with values in $L^2$ (Theorem~\ref{prop:global-kdv}).
Well-posedness in $H^3$ is classical
\citep{BonaSmith1975,Kato1983}; what is new is that it follows from
the scheme itself, together with quantitative error bounds. We do not
address low regularity here.

\paragraph{Organization of the paper.}
The main results are stated in Section~\ref{sec:MainTheo}. The
centered-difference Chernoff function for the Airy group is introduced
in Theorem~\ref{thm:airy-product}. The skew-symmetric splitting yields the
 $L^2$-energy conservation
(Proposition~\ref{prop:mass}), followed by the one-sided $L^2$ stability and uniqueness
(Lemma~\ref{lem:commutator}, \ref{lem:stability} and Corollary~\ref{cor:uniqueness}) and the uniform $H^3$ estimates
(Lemma~\ref{lem:H3} and Proposition~\ref{prop:H3apriori}). The explicit
construction by mesh regularization and Picard iteration is carried
out in Lemma~\ref{lem:positive-mesh} and
Theorem~\ref{thm:wellposed}. The well-posedness is achieved with the
Cauchy property of the approximations (Theorem~\ref{thm:order2}),
uniqueness (Lemma~\ref{lem:kdv-H2-uniqueness}), strong $H^3$
continuity (Proposition~\ref{prop:restart-H3}), the rates of
convergence (Theorem~\ref{thm:kdv-sobolev-quantitative}) and the
proofs of Theorems~\ref{prop:explicit-formula}
and~\ref{prop:global-kdv}. Appendix~\ref{app:interpolation} collects
the interpolation inequalities.

\section{Main result} \label{sec:MainTheo}
{
Let $(T_hf)(x)=f(x+h)$ and, for $h>0$, set
\begin{equation}\label{eq:delta}
 \delta_h=\frac{T_h-T_{-h}}{2h},\qquad \delta_0=D.
\end{equation}
Define the adapted centered Airy generator and its two-parameter unitary propagator by
\begin{equation}\label{eq:Q}
 Q_h=-\delta_h^3+\frac{h^2}{2}\delta_h^5,\quad
 B_n(t)=tQ_{t/n},\quad B_n(0)=0,
\end{equation}
\begin{equation}\label{eq:U}
 U_n(t,s)=\exp(B_n(t)-B_n(s)),\qquad E_n(t)=B_n'(t).
\end{equation}
The nonlinear splitting is\footnote{
{
The identity
$\mathcal F_h(q)=2(T_hq+q+T_{-h}q)\delta_hq$
uses a centered three-point average, of the type encountered in early
finite-difference treatments of KdV \cite{ZabuskyKruskal1965}.
}

}
\begin{equation}\label{eq:Fh}
 \F_h(f)=2f\delta_hf+2\delta_h(f^2).
\end{equation}
The $n$th approximation solves $q_{n,t}=E_n(t)q_n+\F_{t/n}(q_n)$,
$q_n(0)=f$; its positive-mesh regularization is denoted
$p_{n,\varepsilon}$, with $h_{n,\varepsilon}(t)=\max\{t,\varepsilon\}/n$.
}

\begin{theorem}[Constructive representation of the local KdV flow]
\label{prop:explicit-formula}
{
Let\footnote{For $f=0$ the solution is identically $0$ on $\R$.} $0\ne f\in H^3(\R;\R)$, set $T=T_*(f)$, and fix
$0<\varepsilon\le T$, where
\[
T_*(f)=\frac{c_*}{\sqrt{\|f\|_2\|f\|_{H^3}}}
\]
for the universal constant $c_*>0$.
}

%
{
For every integer $n\ge1$, choose a finite subdivision
$0=t_0<\cdots<t_{J(n,\varepsilon)}=T$ 
such that each $I_j=[t_j,t_{j+1}]$ satisfies
\[
\max_j(t_{j+1}-t_j)
\le \frac{\varepsilon}{128A_3n\|f\|_{H^3}},
\]
for some universal constant $A_3>0$.
Set $v_0=f$ and define
}

\begin{equation}\label{eq:Picard-constructive}
 (\Phi_{n,\varepsilon,j}^{v_j}z)(t)
 :=U_n(t,t_j)v_j+\int_{t_j}^tU_n(t,s)
 \F_{h_{n,\varepsilon}(s)}(z(s))\,\dd s,
 \quad z\in C(I_j;H^3),
\end{equation}
where 
$$
h_{n,\varepsilon}(s)=\frac{\max\{s,\varepsilon\}}{n},\qquad 
p_{n,\varepsilon}|_{I_j}:=\lim_{k\to\infty}
 (\Phi_{n,\varepsilon,j}^{v_j})^k(0),\quad \mbox{and}\quad
v_{j+1}:=p_{n,\varepsilon}(t_{j+1}).
$$
Then, \eqref{eq:KdV.init} admits a unique solution 
\begin{equation}\label{eq:KdV-Main}
{
S(\cdot)f\in C([0,T];H^3(\R))
\cap C^1([0,T];L^2(\R)),\qquad}
\end{equation}
given by
\begin{equation}\label{eq:explicit-nonlinear-flow}\multilabel{eq:fixed-epsilon-H3}
 {\displaystyle
 S(\cdot)f=\lim_{n\to\infty} p_{n,\varepsilon}
 \quad\hbox{in }C([0,T];H^2).}
\end{equation}
{and, for every $n\ge1$, $0<\varepsilon\le T$ and
$t\in[0,T]$,}
\begin{equation}\label{eq:explicit-flow-mass}
 \|p_{n,\varepsilon}(t)\|_2
 =\|S(t)f\|_2=\|f\|_2.
\end{equation}

Moreover, for every $f\in H^3(\mathbb R;\mathbb R)$,
\begin{equation}\label{eq:kdv-quant-L2-0}
{
\|p_{n,\varepsilon}-S(\cdot)f\|_{C([0,T];L^2)}
\le C_{T,f}n^{-6/5}.
}
\end{equation}
%
If additionally $f\in H^s$ with $3\leq s\leq 7$, then
\begin{equation}\label{eq:kdv-quant-H2-0}
{
\|p_{n,\varepsilon}-S(\cdot)f\|_{C([0,T];H^2)}
\le C_{T,f,s}n^{-\frac{2(s-2)}5}.
}
\end{equation}


\end{theorem}
\begin{remark}
Taking $\varepsilon=T$ gives the constant nonlinear mesh
$h_{n,T}(t)=T/n$ on $[0,T]$. Formula~\eqref{eq:Picard-constructive}
produces $p_{n,\varepsilon}$ by convergent Picard iterations,
centered translations and explicit Fourier multipliers.
Autonomy and uniqueness of KdV imply $S(t+s)f=S(t)S(s)f$
whenever both sides are defined. 
{The approximate as well as the nonlinear KdV flow preserves the $L^2$-energy.
}
\end{remark}

\begin{theorem}[Global continuation of the constructed KdV flow]
\label{prop:global-kdv}
For any $f\in H^3(\R)$, the local KdV solution constructed in
Theorem~\ref{prop:explicit-formula} extends uniquely to a global solution
$S(\cdot)f\in C(\R;H^3)\cap C^1(\R;L^2)${\color{red}.} 
{
For each $A,\rho>0$ there is $C_{A,\rho}>0$ such that,
whenever $|t|,|s|\le A$ and
$\|f\|_{H^3},\|g\|_{H^3}\le\rho$,
\begin{equation}\label{eq:lipMain}
\|S(t)f-S(t)g\|_{H^2}
\le C_{A,\rho}\|f-g\|_{H^2},
\end{equation}
\begin{equation}\label{eq:global-time-holder}
\|S(t)f-S(s)g\|_{H^2}
\le C_{A,\rho}\bigl(\|f-g\|_{H^2}+|t-s|^{1/3}\bigr),
\end{equation}
and
\begin{equation}\label{eq:global-joint-lipschitz}
\|S(t)f-S(s)g\|_{L^2}
\le C_{A,\rho}\bigl(|t-s|+\|f-g\|_{H^2}\bigr).
\end{equation}
Thus $(t,f)\mapsto S(t)f$ is jointly locally Lipschitz
from $\R\times H^3$ into $L^2$.
}

\end{theorem}

\begin{remark}[Higher regularity and time derivatives]
The high-Sobolev propagation argument in Step~1 of
Theorem~\ref{thm:order2}, combined with the same restart-and-reversibility
argument used for $H^3$, would yield for each integer $s\ge3$ and $f\in H^s$
the local persistence $S(\cdot)f\in C_tH^s$ and 
$S(\cdot)f\in C_t^1H^{s-3}$.
If $s\ge7$, the quadratic $C_tH^2$ rate
is proved in Theorem~\ref{thm:order2}.
\end{remark}

\begin{remark}[Higher-order Richardson--Chernoff corrections]
The present paper establishes the nonlinear second-order convergence rate. However, the (super-)fast Richardson--Chernoff approximation technics for linear evolutions established in \cite{Koudjinan2026SuperfastChernoff} could be used to reach higher arbitrarily prescribed orders for sufficiently regular data.
The construction and proof of those (super-)fast convergent schemes are reserved for a separate work.
\end{remark}
\section{The adapted centered-difference linear evolution}

Let $D=\D$ and consider real-valued solutions of
\begin{equation}\label{eq:kdv}
 u_t=Lu+6uDu,\qquad L=-D^3,\qquad u(0)=u_0,
\end{equation}
on $x\in\R$. For integer $k\ge0$ we use the Sobolev energy norm
$\|f\|_{H^k}^2=\sum_{j=0}^k\|D^jf\|_2^2$.
{
For $r\in\R$ we use the Bessel-potential spaces: $H^r(\R)$ is equipped with the norm
\[
\|f\|_{H^r}^2=\int_\R(1+|\xi|^2)^r|\widehat f(\xi)|^2\,\dd\xi<\infty,
\qquad
\widehat f(\xi)=\frac1{\sqrt{2\pi}}\int_\R e^{-ix\xi}f(x)\,\dd x,
\]
At nonnegative integer
orders, this norm is equivalent to the derivative-sum norm above.
}

{
For real $G,v\in H^k$, the notation set 
$$
\langle G(t),v(t)\rangle_{H^k}:=\sum_{j=0}^k\int_\R D^jG(t,x)\,D^jv(t,x)\,\dd x.$$
 } 
%
Each\footnote{Recall the notations in \S\ref{sec:MainTheo}.} $B_n(t)$ is bounded, skew-adjoint, and these operators commute.\footnote{In fact, they are Fourier multipliers.} Consequently, $U_n(t,\tau)$ is a strongly continuous, unitary evolution family on every $H^r(\R)$:
\begin{equation}\label{eq:unitary}
 \norm[H^r]{U_n(t,\tau)f}=\norm[H^r]{f},\qquad
 U_n(t,\tau)U_n(\tau,\sigma)=U_n(t,\sigma).
\end{equation}

\begin{theorem}[Centered-difference Chernoff function for the Airy group]
\label{thm:airy-product}
Let $L=-D^3$ on $H^r(\R)$, $r\in\R$, and let
$\mathcal S_2(h)=\exp(hQ_h)$ for $h>0$, with $\mathcal S_2(0)=I$.
Then, for every $f\in H^r(\R)$ and $T>0$,
\begin{equation}\label{eq:Airy-Chernoff-product}
 e^{tL}f=\lim_{N\to\infty}\mathcal S_2(t/N)^Nf
 =\lim_{N\to\infty}e^{tQ_{t/N}}f,
 \qquad 0\le t\le T,
\end{equation}
with strong convergence in $H^r$. 
\end{theorem}
\begin{proof}
The Fourier multiplier of $\mathcal S_2(h)^N$ for $h=t/N$ is
$\exp(it h^{-3}g(h\xi))$, where $g(z)=\sin^3z+\tfrac12\sin^5z$.
Since $h^{-3}g(h\xi)\to\xi^3$ as $N\to\infty$, uniformly in $t\in[0,T]$ for each fixed
$\xi$. 
{Since $\widehat{e^{tL}f}(\xi)=e^{it\xi^3}\widehat f(\xi)$} and all multipliers have modulus one, then the dominated convergence and Plancherel's theorem yields the assertion. We refer to \cite{Koudjinan2026SuperfastChernoff} for more informations.
\end{proof}
%
%
For later use, the Fourier symbol of $Q_h$ is $ih^{-3}g(h\xi)$, where
\[
 g(z)=\sin^3z+\frac12\sin^5z
 =z^3-\frac{37}{120}z^7+O(z^9).
\]
Writing $E_n(t)=B_n'(t)$, its symbol for $t>0$ and $h=t/n$ is
\begin{equation}\label{eq:Esymbol}
 \widehat{E_n(t)f}(\xi)
 =\frac{i}{h^3}\bigl(-2g(h\xi)+(h\xi)g'(h\xi)\bigr)\widehat f(\xi).
\end{equation}
At $t=0$, $E_n(0)=L$ on its natural domain. The  global elementary symbol estimates
\[
 |-2g(z)+zg'(z)|\le \frac{19}{2}|z|^3,
 \quad |-2g(z)+zg'(z)-z^3|\le \frac{19}{2}|z|^5,
 \quad |-2g(z)+zg'(z)-z^3|\le \frac{19}{2}|z|^7,
\]
for every $z\in\R$ imply, for $0\le t\le T$ and any $r\in\R$,
\begin{align}
 \norm[H^r]{E_n(t)f}&\le \frac{19}{2}\norm[H^{r+3}]{f},\label{eq:Ebound}\\
 \|(E_n(t)-L)f\|_{H^r}&\le \frac{19}{2}\frac{t^2}{n^2}\|f\|_{H^{r+5}}, \label{eq:Econs2}\\
 \norm[H^r]{(E_n(t)-L)f}
 &\le \frac{19}{2}\frac{t^4}{n^4}\norm[H^{r+7}]{f}.\label{eq:Econs}
\end{align}
\warning{
}
In addition,
\begin{lemma}[Energy identity in the interaction picture]\label{lem:mildenergy}
Let $k\in\{0,1,2,3\}$, $f_0\in H^k$ and
$G\in L^1([0,T];H^k)$. If
\[
 v(t)=U_n(t,0)f_0+\int_0^tU_n(t,\tau)G(\tau)\,\dd\tau,
\]
then $t\mapsto\|v(t)\|_{H^k}^2$ is absolutely continuous and
\begin{equation}\label{eq:mildenergy}
 \frac12\frac{\dd}{\dd t}\|v(t)\|_{H^k}^2
 =\langle G(t),v(t)\rangle_{H^k}\quad\text{for a.e. }t\in(0,T).
\end{equation}
\end{lemma}
\begin{proof}
Set $z(t)=U_n(0,t)v(t)$. The evolution identity and unitarity give
\[
 z(t)=f_0+\int_0^t U_n(0,\tau)G(\tau)\,\dd\tau
 \in W^{1,1}(0,T;H^k).
\]
Consequently $\frac12\frac{\dd}{\dd t}\|z(t)\|_{H^k}^2
=\langle U_n(0,t)G(t),z(t)\rangle_{H^k}$ a.e.
Thus, unitarity of $U_n(t,\tau)$ and its commutation with all $D^j$ imply
\eqref{eq:mildenergy}. This avoids formally differentiating the
unbounded family $E_n(t)$ on $H^k$.
\end{proof}

\section{Skew-symmetric splitting and exact $L^2$-energy conservation}


{
{
Observe that $\F_h(q)=2(T_hq+q+T_{-h}q)\delta_hq$ and is \emph{not} equal to $6q\delta_hq$; rather, it is a skew-symmetric second-order discretization of $6qDq$.
Throughout, ``energy-preserving'' refers to conservation of the $L^2$ norm (Proposition~\ref{prop:mass}); the scheme is not designed to conserve the KdV Hamiltonian $\mathcal{H}(u):=\int(\frac12u_x^2+u^3)\,\dd x$. 
}
}
If $\M_q$ denotes multiplication by a real-valued function $q$, then
\begin{equation}\label{eq:K}
 \F_h(q)=\mathcal K_h(q)q,\qquad
 \mathcal K_h(q)=2\bigl(\M_q\delta_h+\delta_h\M_q\bigr),
 \qquad \mathcal K_h(q)^*=-\mathcal K_h(q).
\end{equation}
{
For the sake of simplicity, we shall write $q_{n,t}:=\partial_tq_n$ for the strong $L^2$ time derivative, i.e., $\lim_{h\to0}\|(q_n(t+h)-q_n(t))/h-q_{n,t}(t)\|_2=0$.\footnote{Of course, with one-sided at endpoints.}} The associated equation is
\begin{equation}\label{eq:scheme}
 q_{n,t}=E_n(t)q_n+\F_{t/n}(q_n),\qquad q_n(0)=u_0\in H^3,
\end{equation}
or, equivalently, the nonlinear integral equation
\begin{equation}\label{eq:mild}
 q_n(t)=U_n(t,0)u_0+
 \int_0^t U_n(t,\tau)\F_{\tau/n}(q_n(\tau))\,\dd\tau.
\end{equation}
We shall use the shorthand\footnote{In particular, $C_tL^2=C([0,T_*];L^2(\R))$,
equipped with the norm
$
\|q\|_{C_tL^2}
:=\sup_{0\le t\le T_*}\|q(t)\|_{L^2(\R)}.
$.}
\[
C_tH^s:=C([0,T_*];H^s(\R)),
\qquad
\|q\|_{C_tH^s}
:=\sup_{0\le t\le T_*}\|q(t)\|_{H^s},
\]
and
\[
L_t^\infty H^s
:=L^\infty([0,T_*];H^s(\R)),
\qquad
\|q\|_{L_t^\infty H^s}
:=
\operatorname*{ess\,sup}_{t\in[0,T_*]}
\|q(t)\|_{H^s(\R)}.
\]

{\color{black}

{
Here and henceforth \emph{strong $C_tL^2$ convergence} of $q_{n,\varepsilon}$ to $q_n$ means
$\lim_{\varepsilon\downarrow0}\sup_{0\le t\le T}\|q_{n,\varepsilon}(t)-q_n(t)\|_{L^2}=0$.
In particular, it is convergence uniformly in time in the spatial $L^2$ norm; it does not assert convergence in $C_tH^3$.
}

{
The equivalence below is the standard variation-of-constants principle for a strongly continuous evolution family;
the proof records the Fourier-multiplier domain justification needed for our time-dependent propagator.}
{
}

\begin{proposition}[Equivalence of the differential and mild formulations]
\label{prop:differential-mild}
Let $n\ge1$, $T>0$, and $u_0\in H^3(\R;\R)$.
Suppose that
\[
q_n\in C([0,T];H^3(\R)).
\]
Then the following assertions are equivalent:
\begin{itemize}
\item[\textbf{(i)}] $q_n$ is a strong $L^2$-valued solution of
\eqref{eq:scheme}, that is,
\[
q_n\in C^1([0,T];L^2(\R)),
\qquad
q_{n,t}=E_n(t)q_n+\F_{t/n}(q_n),
\quad q_n(0)=u_0.
\]

\item[\textbf{(ii)}] $q_n$ satisfies the integral equation
\eqref{eq:mild} in $L^2(\R)$.
\end{itemize}
In particular, every $H^3$-continuous mild solution
belongs to $C^1([0,T];L^2(\R))$.
\end{proposition}

\begin{proof}
We first justify the differentiability of the
evolution family.

Recall that $E_n(t)=B_n'(t)$ and that all
$B_n(t)$ and $E_n(t)$ commute, being Fourier
multipliers. Moreover, by \eqref{eq:Ebound},\footnote{We denote by $C>0$ a universal constant independent of the $t$, $n$ and the initial data.}
\[
\|E_n(t)f\|_{L^2}
\le C\|f\|_{H^3},
\qquad 0\le t\le T,
\]
uniformly in $n$.

The Fourier symbols of $E_n(t)$ are continuous
in $t$, including at $t=0$. 
Thus, by dominated convergence theorem, for every
$f\in H^3(\R)$,
\[
t\longmapsto E_n(t)f
\quad\text{is continuous in }L^2(\R).
\]

It follows, again by Fourier multiplier
differentiation and dominated convergence,
that for every $f\in H^3(\R)$,
\begin{align}
\partial_tU_n(t,\tau)f
&=E_n(t)U_n(t,\tau)f,
\label{eq:U-derivative-t}\\
\partial_\tau U_n(t,\tau)f
&=-U_n(t,\tau)E_n(\tau)f,
\label{eq:U-derivative-tau}
\end{align}
with derivatives taken in $L^2(\R)$.

\smallskip
\noindent
\textbf{(i) $\Longrightarrow$ (ii).}
Assume that $q_n$ solves \eqref{eq:scheme}.
Introduce the interaction variable
\[
v_n(t)=U_n(0,t)q_n(t)
      =e^{-B_n(t)}q_n(t).
\]
Using \eqref{eq:U-derivative-tau},
the product rule in $L^2$, and the
commutativity of the Fourier multipliers,
we obtain
\begin{align*}
v_{n,t}(t)
&=-U_n(0,t)E_n(t)q_n(t)
  +U_n(0,t)q_{n,t}(t)\\
&=U_n(0,t)
  \bigl(q_{n,t}(t)-E_n(t)q_n(t)\bigr)\\
&=U_n(0,t)\F_{t/n}(q_n(t)).
\end{align*}
Since $B_n(0)=0$, we have $v_n(0)=u_0$.
Integrating in $L^2(\R)$ yields
\[
v_n(t)=u_0+
\int_0^t U_n(0,\tau)
\F_{\tau/n}(q_n(\tau))\,\dd\tau.
\]
Multiplying by $U_n(t,0)$ and using
the evolution identity
\[
U_n(t,0)U_n(0,\tau)=U_n(t,\tau),
\]
we obtain
\[
q_n(t)=U_n(t,0)u_0+
\int_0^t U_n(t,\tau)
\F_{\tau/n}(q_n(\tau))\,\dd\tau,
\]
which is precisely \eqref{eq:mild}.

\smallskip
\noindent
\textbf{(ii) $\Longrightarrow$ (i).}
Conversely, assume that $q_n$ satisfies
\eqref{eq:mild} and belongs to
$C([0,T];H^3(\R))$.

Since
\[
\|\delta_hf\|_{H^2}\le\|f\|_{H^3},
\qquad h\in\R,
\]
and $H^3(\R)$ is a Banach algebra,
by the definition \eqref{eq:Fh}, it follows
\[
t\longmapsto\F_{t/n}(q_n(t))
\quad\text{belongs to }
C([0,T];H^2(\R)).
\]
Here we have used the strong continuity of
$h\mapsto\delta_h$ from $H^3$ into $H^2$,
including at $h=0$.

Multiplying \eqref{eq:mild} by $U_n(0,t)$
yields
\[
v_n(t):=U_n(0,t)q_n(t)
=u_0+\int_0^t U_n(0,\tau)
\F_{\tau/n}(q_n(\tau))\,\dd\tau.
\]
The integrand is continuous in $H^2$.
Hence the fundamental theorem of calculus
for Bochner integrals implies
\[
v_n\in C^1([0,T];H^2(\R))
\]
and
\[
v_{n,t}(t)
=U_n(0,t)\F_{t/n}(q_n(t)).
\]

On the other hand,
\[
q_n(t)=U_n(t,0)v_n(t),
\qquad
v_n\in C([0,T];H^3(\R)).
\]
Using \eqref{eq:U-derivative-t} and the
product rule in $L^2$, we conclude that
$q_n\in C^1([0,T];L^2(\R))$ and
\begin{align*}
q_{n,t}(t)
&=E_n(t)U_n(t,0)v_n(t)
  +U_n(t,0)v_{n,t}(t)\\
&=E_n(t)q_n(t)
  +U_n(t,0)U_n(0,t)
       \F_{t/n}(q_n(t))\\
&=E_n(t)q_n(t)+\F_{t/n}(q_n(t)).
\end{align*}
Finally, evaluating \eqref{eq:mild} at
$t=0$ yields $q_n(0)=u_0$.
Thus $q_n$ satisfies \eqref{eq:scheme},
which completes the proof.
\end{proof}
}

\begin{proposition}[Exact $L^2$ conservation]\label{prop:mass}
{
Fix $n\ge1$ and $T>0$. Every real solution $q\in C([0,T];H^3(\R))\cap C^1([0,T];L^2(\R))$ of \eqref{eq:scheme}, with the differential equation understood in $L^2$ and $q(0)=u_0\in H^3$, satisfies
\begin{equation}\label{eq:mass}
 \norm[L^2]{q(t)}=\norm[L^2]{u_0},\qquad 0\le t\le T.
\end{equation}

}
\end{proposition}
\begin{proof}
{
{
For real $q$, 
$\mathcal{F}_h(q)=\mathcal K_h(q)q$ with
$\mathcal K_h(q)=2(\M_q\delta_h+\delta_h\M_q)$, which} 
is skew-adjoint, by \eqref{eq:K}.
The same holds for $E_n(t)$ 
{on $H^3$}. 
Hence the assumed $C^1_tL^2\cap C_tH^3$ regularity justifies
$\frac12\frac{\dd}{\dd t}\|q(t)\|_2^2=\langle E_n(t)q(t),q(t)\rangle_2+\langle 2(M_q\mathcal A+\mathcal A M_q)q,q\rangle_2=0$.
}
\end{proof}
Next, we compute the Hamiltonian drift for the approximants.
{
\begin{proposition}[Hamiltonian drift for the approximations]
\label{prop:discrete-H-balance}
Let $\mathcal H(q)=\int_\R(\tfrac12|Dq|^2+q^3)\,\dd x$.
Suppose $q\in C([0,T];H^3)\cap C^1([0,T];L^2)$ is real,
$q_t=E_n(t)q+\mathcal F_{h(t)}(q)$ in $L^2$, and
$h:[0,T]\to[0,\infty)$ is continuous, with $\delta_0=D$.
Then
\begin{equation}\label{eq:approx-H-balance}
\frac{\dd}{\dd t}\mathcal H(q(t))
=\left\langle-D^2q+3q^2,
(E_n(t)-L)q+\mathcal F_{h(t)}(q)-6qDq\right\rangle_{L^2}.
\end{equation}
If $\sup_{[0,T]}\|q(t)\|_{H^5}\le M_5$, then the Hamiltonian drift is:
\begin{equation}\label{eq:approx-H-defect}
|\mathcal H(q(t))-\mathcal H(q(0))|
\le C_{M_5,T}\int_0^t
\left(\frac{s^2}{n^2}+h(s)^2\right)\,\dd s.
\end{equation}
\end{proposition}
\begin{proof}
Since $\mathcal H'(q)=-D^2q+3q^2\in L^2$,
$Lq+6qDq=D\mathcal H'(q)$, and
\begin{equation}\label{eq:KdV-fiels}
\langle\mathcal H'(q),D\mathcal H'(q)\rangle_2=0,
\end{equation}
subtracting this term yields the drift identity.
By \eqref{eq:Econs2},
$\|(E_n(s)-L)q\|_2\le Cs^2n^{-2}\|q\|_{H^5}$;
Taylor's formula for centered differences gives
$\|\mathcal F_h(q)-6qDq\|_2\le Ch^2\|q\|_{H^3}^2$.
The $L^2$ norm of $\mathcal H'(q)$ is bounded by a polynomial in
$M_5$, so integration proves \eqref{eq:approx-H-defect}.
\end{proof}
}


{
\section{One-sided $L^2$ stability without derivative loss}
}

\begin{lemma}[Discrete commutator]\label{lem:commutator}
For  $a\in W^{1,\infty}(\R,\R)$, $f\in L^2$, and $h\in\R$,
\begin{equation}\label{eq:commutator}
 \norm[L^2]{[a,\delta_h]f}\le \norm[L^\infty]{Da}\norm[L^2]{f}, 
 \qquad
 \bigl|\ip{a\delta_hf}{f}_{L^2}\bigr|
 \le \frac12\norm[L^\infty]{Da}\norm[L^2]{f}^2.
\end{equation}
\end{lemma}
\begin{proof}
For $h\ne0$,
\[
 [a,\delta_h]f(x)
 =\frac{a(x)-a(x+h)}{2h}f(x+h)
 -\frac{a(x)-a(x-h)}{2h}f(x-h).
\]
Translation invariance of $L^2$ and the mean-value theorem give the first estimate. For the second estimate, from $\delta_h^*=-\delta_h$ it follows
$2\ip{a\delta_hf}{f}=\ip{[a,\delta_h]f}{f}$ on $L^2$. Then, letting $h\to0$ complete the proof.
\end{proof}

\begin{lemma}[{Polarization and one-sided Lipschitz estimate}]
\label{lem:stability}
For $q,p\in H^2(\R,\R)$ and $w=q-p$, one has
\begin{equation}\label{eq:polar}
 \F_h(q)-\F_h(p)
 =(q+p)\delta_hw+w\delta_h(q+p)+2\delta_h\bigl((q+p)w\bigr).
\end{equation}
In particular, for every $h\in\R$,
\begin{equation}\label{eq:onesided}
 \bigl|\ip{\F_h(q)-\F_h(p)}{q-p}_{L^2}\bigr|
 \le \frac32\norm[L^\infty]{D(q+p)}\norm[L^2]{q-p}^2.
\end{equation}
\end{lemma}
\begin{proof}
Set $a=q+p$. Using the skew-adjointness of $\delta_h$ in \eqref{eq:polar},
\[
 \ip{\F_h(q)-\F_h(p)}{w}
 =-\ip{a\delta_hw}{w}+\ip{(\delta_ha)w}{w},
\]
which, together with \eqref{eq:commutator} and
$\norm[L^\infty]{\delta_ha}\le\norm[L^\infty]{Da}$
prove \eqref{eq:onesided}.
\end{proof}

\begin{corollary}[Uniqueness and stability]\label{cor:uniqueness}
{
Let $n\ge1$. Suppose $q,p\in C([0,T];L^2(\R))\cap L^\infty([0,T];H^3(\R))$ are real-valued mild solutions of \eqref{eq:mild}, possibly with different initial data. 
Then
\begin{equation}\label{eq:stability}
 \norm[L^2]{q(t)-p(t)}
 \le\norm[L^2]{q(0)-p(0)}
 \exp\left(\frac32\int_0^t\norm[L^\infty]{D(q(\tau)+p(\tau))}\,\dd\tau\right).
\end{equation}
In particular, there is at most one mild solution in this class with prescribed initial data.
}
\end{corollary}
\begin{proof}
Let $q,p\in L^\infty([0,T];H^3)\cap C([0,T];L^2)$ be mild solutions of \eqref{eq:mild}. By interpolation, $q,p\in C([0,T];H^\sigma)$ for every $\sigma<3$; since $h\mapsto\delta_h$ is strongly continuous from $H^1$ to $L^2$, the forcing $G(\tau)=\F_{\tau/n}(q(\tau))-\F_{\tau/n}(p(\tau))$ belongs to $C([0,T];L^2)$. The difference $w=q-p$ solves the mild equation with initial value $q(0)-p(0)$ and forcing $G$. Lemma~\ref{lem:mildenergy} with $k=0$ and \eqref{eq:onesided} give
\[
 \tfrac12\tfrac{\dd}{\dd t}\norm[L^2]{w}^2\le\tfrac32\norm[L^\infty]{D(q+p)}\norm[L^2]{w}^2
 \quad\text{a.e.},
\]
and Gronwall's inequality yields \eqref{eq:stability}.
\end{proof}

\section{A uniform $H^3$ estimate enhanced by $L^2$-energy conservation}

\begin{lemma}[Uniform $H^3$ nonlinear energy estimate]\label{lem:H3}
For real $q\in H^3(\R)$ and every $h\in\R\setminus\{0\}$,
\begin{equation}\label{eq:H3nonlinear}
 \left|\sum_{k=0}^3\ip{D^k\F_h(q)}{D^kq}_{L^2}\right|
 \le C\norm[L^\infty]{Dq}\norm[H^3]{q}^2,
\end{equation}
where $C$ is independent of $h$. The same estimate holds at $h=0$ for smooth $q$, but the pairing at $k=3$ is not necessarily defined on all of $H^3$.
\end{lemma}
\begin{proof}
For $k=0$, the pairing vanishes by skew-adjointness. For $k\ge1$, set $f=D^kq$ and expand
\begin{align*}
 D^k\F_h(q)
 ={}&2q\delta_hf+4\delta_h(qf)
 +2\sum_{j=1}^k\binom kj(D^jq)\delta_hD^{k-j}q\\
 &+2\delta_h\left(\sum_{j=1}^{k-1}\binom kj(D^jq)(D^{k-j}q)\right).
\end{align*}
The first two terms satisfy
\[
 \ip{2q\delta_hf+4\delta_h(qf)}{f}
 =-2\ip{q\delta_hf}{f},
\]
whose absolute value is bounded by $\norm[L^\infty]{Dq}\norm[L^2]{f}^2$ by Lemma~\ref{lem:commutator}.
We also have the mesh-uniform estimates
\[
 \norm[L^p]{\delta_hg}\le\norm[L^p]{Dg},\qquad 1\le p\le\infty.
\]
For $k=1,2$ all remaining products are bounded by
$C\norm[L^\infty]{Dq}\norm[H^3]{q}^2$ using 
{H\"older's} inequality. For $k=3$, the middle products additionally require
\begin{equation}\label{eq:interpolation-special}
 \norm[L^4]{D^2q}^2\le C\norm[L^\infty]{Dq}\norm[L^2]{D^3q}.
\end{equation}
To see this, put $f=Dq$, integrate $\int |f'|^4$ by parts for smooth compactly supported $f$, and estimate
\[
 \norm[L^4]{f'}^4
 \le3\norm[L^\infty]{f}\norm[L^4]{f'}^2\norm[L^2]{f''}.
\]
Approximation proves \eqref{eq:interpolation-special} in $H^3$.
Thus, 
\[
 \norm[L^2]{(D^2q)\delta_h(Dq)}
 \le\norm[L^4]{D^2q}\norm[L^4]{\delta_hDq}
 \le C\norm[L^\infty]{Dq}\norm[L^2]{D^3q},
\]
and
\[
 \norm[L^2]{\delta_h((Dq)(D^2q))}
 \le\norm[L^2]{D((Dq)(D^2q))}
 \le C\norm[L^\infty]{Dq}\norm[L^2]{D^3q}.
\]
These estimates cover all the terms at $k=3$ and prove the lemma.
\end{proof}

\begin{proposition}[Uniform a priori bound]\label{prop:H3apriori}
{
Let $h:[0,T]\to(0,\infty)$ be continuous, and let $q\in C([0,T];H^3(\R;\R))$ be a mild solution (equivalently, an $L^2$-strong solution by Proposition~\ref{prop:differential-mild}, with the same proof for $h$) of}
\[
 q_t=E_n(t)q+\F_{h(t)}(q),\qquad q(0)=u_0\in H^3,
\]
Assume\footnote{For $u_0=0$ the solution is identically zero.} $u_0\ne0$ and 
let $M_0=\norm[L^2]{u_0}$ and $Y_0=\norm[H^3]{u_0}$. There exists an absolute $C>0$, independent of $n$, $t$ and the mesh $h$, such that
\begin{equation}\label{eq:Ydiff}
 Y'(t)\le C M_0^{1/2}Y(t)^{3/2},
 \qquad Y(t)=\norm[H^3]{q(t)},
\end{equation}
and hence
\begin{equation}\label{eq:Ybound}\multilabel{eq:Tstar}
 Y(t)\le\frac{Y_0}{\left(1-\frac C2 t\sqrt{M_0Y_0}\right)^2}, \qquad {0\leq t\leq\min\{T,T_*\}},\qquad T_* :=\frac{1}{C\sqrt{M_0Y_0}},
\end{equation}
 so that 
$\sup_{t\in[0,\min\{T,T_*\}]}\norm[H^3]{q(t)}\le 4\norm[H^3]{u_0}$.
\end{proposition}
\begin{proof}
The linear evolution is unitary on $H^3$ and its generator is skew-adjoint. Equivalently, write the energy in the interaction variable $U_n(0,t)q(t)$, thereby avoiding the unbounded-generator domain issue. Lemma~\ref{lem:H3} gives
\[
 \frac12\frac{\dd}{\dd t}Y(t)^2
 \le C\norm[L^\infty]{Dq(t)}Y(t)^2.
\]
The one-dimensional Gagliardo--Nirenberg inequality gives
\[
 \norm[L^\infty]{Dq}
 \le C\norm[L^2]{q}^{1/2}\norm[L^2]{D^3q}^{1/2}
 \le C M_0^{1/2}Y^{1/2},
\]
where Proposition~\ref{prop:mass} supplies $L^2$-conservation. Dividing by $Y$ and integrating the resulting scalar differential inequality proves \eqref{eq:Ydiff}--\eqref{eq:Ybound}.
\end{proof}

\section{Explicit construction via temporal mesh regularization}

{
Throughout Sections~6--7 we write
$q_{n,\varepsilon}=p_{n,\varepsilon}$.
}
For $\varepsilon>0$ set
\begin{equation}\label{eq:meshreg}
 h_{n,\varepsilon}(t)=\frac{\max\{t,\varepsilon\}}n, \qquad \Lambda=\frac{n}{\varepsilon}.
\end{equation}
Consider the regularized equation\footnote{Only the \emph{nonlinear} mesh is regularized, $U_n$ is unchanged.}
\begin{equation}\label{eq:reg}
 q_{n,\varepsilon}(t)=U_n(t,0)u_0
 +\int_0^t U_n(t,\tau)
 \F_{h_{n,\varepsilon}(\tau)}(q_{n,\varepsilon}(\tau))\,\dd\tau, \qquad q_{n,\varepsilon}(0)=u_0.
\end{equation}

\begin{lemma}[Contraction for positive mesh]\label{lem:positive-mesh}
For each fixed $n\ge1$, $\varepsilon>0$, and $u_0\in H^3$, equation \eqref{eq:reg} has a unique local solution $q_{n,\varepsilon}$ in $C([0,T_{n,\varepsilon}];H^3)$. Its lifespan can be extended to the common interval $[0,T_*]$ from \eqref{eq:Tstar}, with
\begin{equation}\label{eq:uniformreg}
 \sup_{n\ge1,\varepsilon>0}\sup_{0\le t\le T_*}
 \norm[H^3]{q_{n,\varepsilon}(t)}\le 4\norm[H^3]{u_0}.
\end{equation}
\end{lemma}

\begin{proof}
We divide the proof into three steps.

\medskip
\noindent
\textbf{Step 1. Local Lipschitz continuity.}

Since $h(t)\ge\varepsilon/n$, the Fourier
multiplier representation of $\delta_h$ yields
\begin{equation}\label{eq:positive-mesh-bound}
\|\delta_{h(t)}f\|_{H^3}
\le \Lambda\|f\|_{H^3},
\qquad t\ge0.
\end{equation}
Let $A_3\ge1$ be a constant such that
\[
\|fg\|_{H^3}
\le A_3\|f\|_{H^3}\|g\|_{H^3}.
\]
Recall that
$
\F_h(q)=2q\delta_hq+2\delta_h(q^2).
$ 
Then, for every $q\in H^3$,
\begin{align}
\|\F_{h(t)}(q)\|_{H^3}
&\le
2A_3\|q\|_{H^3}
       \|\delta_{h(t)}q\|_{H^3}
 +2\|\delta_{h(t)}(q^2)\|_{H^3}
\nonumber\\
&\le
4A_3\Lambda\|q\|_{H^3}^2.
\label{eq:positive-mesh-growth}
\end{align}

Furthermore, for $q,p\in H^3$, the identity
\[
\begin{aligned}
\F_h(q)-\F_h(p)
={}&2(q-p)\delta_hq
   +2p\delta_h(q-p)+2\delta_h\bigl((q-p)(q+p)\bigr)
\end{aligned}
\]
implies
\begin{equation}\label{eq:positive-mesh-Lip}
\|\F_{h(t)}(q)-\F_{h(t)}(p)\|_{H^3}
\le
4A_3\Lambda
\bigl(\|q\|_{H^3}+\|p\|_{H^3}\bigr)
\|q-p\|_{H^3}.
\end{equation}

Thus $\F_{h(t)}:H^3\to H^3$ is locally
Lipschitz, uniformly in $t$, with constants
depending on $n/\varepsilon$.

Moreover, $t\mapsto\delta_{h(t)}f$ is
continuous in $H^3$ for each $f\in H^3$.
Indeed, the corresponding Fourier multipliers
converge pointwise as $t$ varies and are
uniformly bounded by $\Lambda$. Dominated
convergence gives strong continuity.
It follows from \eqref{eq:positive-mesh-Lip}
that
\[
(t,q)\longmapsto\F_{h(t)}(q)
\]
is jointly continuous from
$[0,\infty)\times H^3$ into $H^3$.

\medskip
\noindent
\textbf{Step 2. Banach's contraction theorem.}

Let $t_0\ge0$ and $v_0\in H^3(\R;\R)$
be arbitrary. We consider the integral
equation with initial value $v_0$ at $t_0$:
\[
q(t)=U_n(t,t_0)v_0+
\int_{t_0}^t U_n(t,\tau)
\F_{h(\tau)}(q(\tau))\,\dd\tau.
\]

Assume first that $v_0\ne0$ and set
\[
M=\|v_0\|_{H^3},
\qquad R=2M.
\]
For $\Delta>0$, define
\[
X_\Delta=C([t_0,t_0+\Delta];H^3(\R;\R))
\]
and its closed ball
\[
\mathcal B_R
=\left\{q\in X_\Delta:
\|q\|_{X_\Delta}\le R\right\}, \qquad 
\|q\|_{X_\Delta}
:=\sup_{t\in[t_0,t_0+\Delta]}
\|q(t)\|_{H^3(\R)}.
\]
This is a complete metric space.

Define $\Phi:\mathcal B_R\to X_\Delta$ by
\[
(\Phi q)(t)
=
U_n(t,t_0)v_0+
\int_{t_0}^t U_n(t,\tau)
\F_{h(\tau)}(q(\tau))\,\dd\tau.
\]
The strong continuity of $U_n$, together
with Step~1, ensures that $\Phi q$ belongs
to $X_\Delta$.

Since $U_n(t,\tau)$ is unitary on $H^3$,
\eqref{eq:positive-mesh-growth} yields
\begin{align*}
\|\Phi q\|_{X_\Delta}
&\le M+
\int_{t_0}^{t_0+\Delta}
\|\F_{h(\tau)}(q(\tau))\|_{H^3}
\,\dd\tau \le M+4A_3\Lambda\Delta R^2.
\end{align*}

Similarly, for $q,p\in\mathcal B_R$,
\begin{align*}
\|\Phi q-\Phi p\|_{X_\Delta}
&\le
\int_{t_0}^{t_0+\Delta}
\|\F_{h(\tau)}(q(\tau))
-\F_{h(\tau)}(p(\tau))\|_{H^3}
\,\dd\tau\le
8A_3\Lambda R\Delta
\|q-p\|_{X_\Delta}.
\end{align*}

Choose
\begin{equation}\label{eq:positive-mesh-time}
0<\Delta\le
\frac{1}{16A_3\Lambda R}.
\end{equation}
Then, since $M=R/2$,
\[
\|\Phi q\|_{X_\Delta}
\le\frac R2+\frac R4<R,
\]
and
\[
\|\Phi q-\Phi p\|_{X_\Delta}
\le\frac12\|q-p\|_{X_\Delta}.
\]

Thus $\Phi$ maps $\mathcal B_R$ into itself
and is a contraction. Banach's
fixed-point theorem yields a unique
solution in $\mathcal B_R$.

Local uniqueness among arbitrary continuous
$H^3$ solutions follows by applying
\eqref{eq:positive-mesh-Lip} and Gronwall's
inequality on intervals where their
$H^3$ norms are bounded.

If $v_0=0$, the zero function is a solution,
and the same Lipschitz estimate gives
uniqueness.

We have therefore established local
well-posedness in $H^3$ for every initial
time $t_0\ge0$ and every initial datum
$v_0\in H^3$.

In particular, taking $t_0=0$ and
$v_0=u_0$ gives a unique maximal solution
\[
q_{n,\varepsilon}
\in C([0,T_{n,\varepsilon}^{\max});H^3),
\]
where $T_{n,\varepsilon}^{\max}>0$.

\medskip
\noindent
\textbf{Step 3. Uniform a priori bound and continuation.}

Proposition~\ref{prop:mass} implies
\[
\|q_{n,\varepsilon}(t)\|_{L^2}
=\|u_0\|_{L^2}.
\]

Moreover, Proposition~\ref{prop:H3apriori} implies
\begin{equation}\label{eq:positive-mesh-uniform}
\sup_{0\le t<\min\{
T_{n,\varepsilon}^{\max},T_*\}}
\|q_{n,\varepsilon}(t)\|_{H^3}
\le 4\|u_0\|_{H^3},
\end{equation}
where $T_*$ is as in \eqref{eq:Tstar}.

It remains to justify continuation up
to $T_*$. recall $Y_0=\|u_0\|_{H^3}$. 
Suppose, by contradiction,
that
\[
T_{n,\varepsilon}^{\max}< T_*.
\]
By \eqref{eq:positive-mesh-growth}
and \eqref{eq:positive-mesh-uniform},
%
the nonlinear forcing is bounded in
$H^3$ on the maximal interval:
\[
\|\F_{h(t)}(q_{n,\varepsilon}(t))\|_{H^3}
\le {64}A_3\Lambda Y_0^2.
\]
Introduce
\[
z(t)=U_n(0,t)q_{n,\varepsilon}(t).
\]
The evolution identity yields
\[
z(t)=u_0+
\int_0^tU_n(0,\tau)
\F_{h(\tau)}(q_{n,\varepsilon}(\tau))
\,\dd\tau.
\]
%
By unitarity,
\[
\|z(t)-z(s)\|_{H^3}
\le {64}A_3\Lambda Y_0^2|t-s|
\]
for $s,t<T_{n,\varepsilon}^{\max}$.

Hence $z(t)$ has a strong limit in
$H^3$ as
$t\uparrow T_{n,\varepsilon}^{\max}$.
The strong continuity of $U_n(t,0)$
then implies that
$q_{n,\varepsilon}(t)$ also has a
strong $H^3$ limit, say $v_*$.

Applying Step~2 with initial time
$T_{n,\varepsilon}^{\max}$ and datum
$v_*$ extends the solution beyond
$T_{n,\varepsilon}^{\max}$, contradicting
maximality.

Therefore
\[
T_{n,\varepsilon}^{\max}\ge T_*.
\]
Combining this with
\eqref{eq:positive-mesh-uniform} proves
\[
\sup_{n\ge1,\varepsilon>0}
\sup_{0\le t\le T_*}
\|q_{n,\varepsilon}(t)\|_{H^3}
\le 4Y_0,
\]
which is \eqref{eq:uniformreg}.
\end{proof}

\begin{lemma}[Control of the mesh variation]\label{lem:meshvariation}
For all $a,b\ge0$, $f\in H^3(\R)$,
\begin{equation}\label{eq:meshvariation}
 \norm[L^2]{(\delta_a-\delta_b)f}
 \le C|a^2-b^2|\norm[H^3]{f}.
\end{equation}
Consequently, for $p\in H^3(\R;\R)$,
\begin{equation}\label{eq:Fmesh}
 \norm[L^2]{(\F_a-\F_b)(p)}
 \le C|a^2-b^2|\norm[H^3]{p}^2.
\end{equation}
\end{lemma}
\begin{proof}
The Fourier symbol of $\delta_h$ is $d_h(\xi)=i\sin(h\xi)/h$. For $h>0$,
\[
 \partial_hd_h(\xi)
 =i\frac{h\xi\cos(h\xi)-\sin(h\xi)}{h^2},
 \quad |\partial_hd_h(\xi)|\le Ch|\xi|^3.
\]
The last inequality follows from the global bound
$|z\cos z-\sin z|\le C|z|^3$. Integrating between $a$ and $b$ and applying Plancherel gives \eqref{eq:meshvariation}. Then use
\[
 (\F_a-\F_b)(p)
 =2p(\delta_a-\delta_b)p+2(\delta_a-\delta_b)(p^2)
\]
and the algebra property of $H^3$ to obtain \eqref{eq:Fmesh}.
\end{proof}

\begin{proposition}[Quantitative removal of the regularization]\label{prop:Cauchy}
For $0<\varepsilon<\eta\le T_*$ and 
{$R=4Y_0$},
\begin{equation}\label{eq:epsCauchy}
 \sup_{0\le t\le T_*}\norm[L^2]{q_{n,\varepsilon}(t)-q_{n,\eta}(t)}
 \le C_{R,T_*}\frac{\eta^3-{\varepsilon}^3}{n^2}.
\end{equation}
In particular, $\{q_{n,\varepsilon}\}$ converges (as ${\varepsilon\downarrow0}$)  to some $q_n$ in $C([0,T_*];L^2)$, uniformly in $n$.
\end{proposition}
\begin{proof}
Let $w=q_{n,\varepsilon}-q_{n,\eta}$. Subtraction gives
\begin{align*}
 w_t={}&E_n(t)w
 +\F_{h_{n,\varepsilon}}(q_{n,\varepsilon})
 -\F_{h_{n,\varepsilon}}(q_{n,\eta})+\bigl(\F_{h_{n,\varepsilon}}-\F_{h_{n,\eta}}\bigr)(q_{n,\eta}).
\end{align*}
Lemma~\ref{lem:stability} and \eqref{eq:uniformreg}, together with the embedding $H^3\hookrightarrow W^{1,\infty}$, control the first nonlinear difference in the $L^2$ energy. 
{The sole remaining forcing term is controlled by \eqref{eq:Fmesh}; the $L^2$-energy differentiation is justified by Lemma~\ref{lem:mildenergy} with $k=0$. 
} 
Then, 
\[
 \frac{\dd}{\dd t}\norm[L^2]{w(t)}
 \le CR\norm[L^2]{w(t)}
 +CR^2|h_{n,\varepsilon}(t)^2-h_{n,\eta}(t)^2|.
\]
Since $w(0)=0$, Gronwall's inequality
yields
\[
\|w\|_{C_tL^2}
\le CR^2e^{CRT_*}
\int_0^{T_*}
|h_{n,\varepsilon}(t)^2
-h_{n,\eta}(t)^2|\,\dd t.
\]

Observe that the two meshes coincide
for $t\ge\eta$. Hence
\begin{align*}
&\int_0^{T_*}
|h_{n,\varepsilon}(t)^2
-h_{n,\eta}(t)^2|\,\dd t\\
&\quad=
\frac1{n^2}
\left[
\int_0^\varepsilon(\eta^2-\varepsilon^2)\,\dd t
+\int_\varepsilon^\eta(\eta^2-t^2)\,\dd t
\right]\\
&\quad=\frac{2}{3n^2}(\eta^3-\varepsilon^3),
\end{align*}
whence \eqref{eq:epsCauchy}. 
It follows that $\{q_{n,\varepsilon}\}_{\varepsilon\downarrow0}$ is Cauchy in $C([0,T_*];L^2)$. Therefore, $q_n=\lim_{\varepsilon\downarrow0} (q_{n,\varepsilon})\in C([0,T_*];L^2)$, 
{uniformly} in $n$.
\end{proof}

\begin{theorem}[Well-posedness]\label{thm:wellposed}
Let $0\not= u_0\in H^3(\R;\R)$. 
For every integer $n\ge1$, equation \eqref{eq:mild} has a unique real solution
\begin{equation}\label{eq:solutionreg}
 q_n\in C([0,T_*];H^3(\R)),\qquad  T_* :=\frac{1}{C\sqrt{\norm[L^2]{u_0}\norm[H^3]{u_0}}},
\end{equation}
and satisfies, for every $t\in [0,T_*]$,

\begin{align}
 \norm[L^2]{q_n(t)}&=\norm[L^2]{u_0},\label{eq:exactmassfinal}\\
 \sup_{n\ge1}\sup_{0\le t\le T_*}\norm[H^3]{q_n(t)}&\le {4}\norm[H^3]{u_0},\label{eq:H3final}\\
 \sup_{0\le t\le T_*}\norm[L^2]{q_{n,\varepsilon}(t)-q_n(t)}
 &\le C\frac{\varepsilon^3}{n^2},\quad 0<\varepsilon\le T_*.
 \label{eq:reg-rate}
\end{align}
\end{theorem}
\begin{proof}
By Proposition~\ref{prop:Cauchy}, $q_{n,\varepsilon}\to q_n$ in $C_tL^2$, and the uniform $H^3$ bound \eqref{eq:uniformreg} yields $q_n(t)\in H^3$ for every $t$, with the asserted bound. Sobolev interpolation shows that, for each $\sigma<3$,
\begin{equation}\label{eq:interpolation}
 q_{n,\varepsilon}\longrightarrow q_n
 \quad\text{in }C([0,T_*];H^\sigma).
\end{equation}
In particular, this convergence holds in $C_tH^2$. 
{

{Fix $0<\eta\le T_*$ and let the \emph{smaller} parameter
$\varepsilon\downarrow0$ in \eqref{eq:epsCauchy}, using
$q_{n,\varepsilon}\to q_n$ in $C_tL^2$. Renaming the fixed parameter
$\eta$ as $\varepsilon$ gives
\begin{equation}\label{eq:reg-quantitative}
\|q_{n,\varepsilon}-q_n\|_{C_tL^2}
\le C_{R,T_*}\frac{\varepsilon^3}{n^2}.
\end{equation}
Since $\sup_{t}\bigl(\|q_{n,\varepsilon}(t)\|_{H^3}
+\|q_n(t)\|_{H^3}\bigr)\le8Y_0$,
Sobolev interpolation yields
\begin{align}
\|q_{n,\varepsilon}-q_n\|_{C_tH^2}
&\le C\|q_{n,\varepsilon}-q_n\|_{C_tL^2}^{1/3}
       \|q_{n,\varepsilon}-q_n\|_{L_t^\infty H^3}^{2/3}
\le C_{T_*,u_0}\varepsilon n^{-2/3}.
\label{eq:reg-H2-quantitative}
\end{align}
The constants are uniform in $n$ and $0<\varepsilon\le T_*$.}

Now, we prove the convergence of the 
{nonlinear} part. 
Fix $0<t\leq T_*$. For $0<\varepsilon<t$,
\[
h_{n,\varepsilon}(t)=t/n.
\]
Since\footnote{By analogy to \eqref{eq:positive-mesh-Lip}, taking the $H^2$-norm instead of $H^3$.} $\F_{t/n}:H^2\to H^2$ is
locally Lipschitz for fixed $t>0$ and $n\geq 1$,
\eqref{eq:interpolation} implies as $\varepsilon\downarrow0$,
\[
\F_{h_{n,\varepsilon}(t)}
(q_{n,\varepsilon}(t))=\F_{t/n}
(q_{n,\varepsilon}(t))
\longrightarrow
\F_{t/n}(q_n(t))
\quad\text{in }H^2.
\]

Furthermore, uniformly in $h\ge0$,
\[
\|\F_h(f)\|_{H^2}
\le C\|f\|_{H^3}^2.
\]
The uniform $H^3$ bounds therefore provide
an integrable dominating function.

By the 
dominated convergence
theorem, as $\varepsilon\downarrow0$,
\[
\F_{h_{n,\varepsilon}(\cdot)}
(q_{n,\varepsilon})
\longrightarrow
\F_{\cdot/n}(q_n)
\quad\text{in }L^1([0,T_*];H^2).
\]

Finally, the unitarity of $U_n(t,\tau)$
on $H^2$ yields as $\varepsilon\downarrow0$
\begin{align*}
&\sup_{0\le t\le T_*}
\left\|
\int_0^tU_n(t,\tau)
\Bigl[
\F_{h_{n,\varepsilon}(\tau)}
(q_{n,\varepsilon}(\tau))
-\F_{\tau/n}(q_n(\tau))
\Bigr]\,\dd\tau
\right\|_{H^2}\\
&\quad\le
\int_0^{T_*}
\left\|
\F_{h_{n,\varepsilon}(\tau)}
(q_{n,\varepsilon}(\tau))
-\F_{\tau/n}(q_n(\tau))
\right\|_{H^2}\,\dd\tau
\longrightarrow0.
\end{align*}

We may therefore pass to the limit
in the regularized mild equation,
obtaining \eqref{eq:mild} in $H^2$.

}
The strong $L^2$ limit preserves the $L^2$-energy and the quantitative bound \eqref{eq:reg-rate}.

Uniqueness in the class $L^\infty_tH^3\cap C_tL^2$ follows from Corollary~\ref{cor:uniqueness}.

{

We finally prove the strong continuity
$
q_n\in C([0,T_*];H^3(\R)).
$

\noindent
\textbf{Step 1. Strong continuity for positive times.}
Fix $t_0\in(0,T_*)$. Since
$q_n(t_0)\in H^3(\R)$, choose
$0<\rho<t_0/2$ such that
\[
I=[t_0-\rho,t_0+\rho]
\subset(0,T_*).
\]
For every $t\in I$, we have
\[
h(t)=\frac tn\ge\frac{t_0}{2n}>0.
\]
 Consequently,
\[
\norm[H^3]{\delta_{t/n}f}
\le\frac{2n}{t_0}\norm[H^3]{f}.
\]
Consider the integral equation
\begin{equation}\label{eq:local-around-t0}
v(t)=U_n(t,t_0)q_n(t_0)
+\int_{t_0}^t
U_n(t,\tau)\F_{\tau/n}(v(\tau))
\,\dd\tau.
\end{equation}
The same contraction argument as in
Lemma~\ref{lem:positive-mesh}
yields, after possibly reducing $\rho$,
a unique solution
$
v\in C(I;H^3(\R)).
$

On the other hand, the evolution identity
$
U_n(t,\tau)U_n(\tau,s)=U_n(t,s)
$ 
and the mild equation \eqref{eq:mild}
show that $q_n$ also satisfies
\eqref{eq:local-around-t0} as an identity
in $L^2$.

We claim that these two solutions coincide.
Indeed, for $a,b$ in a bounded subset
of $H^3$ and $t\in I$, the identity
\[
\F_h(a)-\F_h(b)
=
2(a-b)\delta_ha
+2b\delta_h(a-b)
+2\delta_h((a+b)(a-b))
\]
together with
\[
\norm[L^2]{\delta_hf}
\le h^{-1}\norm[L^2]{f},
\]
 and the embedding
$H^3\hookrightarrow W^{1,\infty}$, it follows
\[
\norm[L^2]{\F_{t/n}(a)-\F_{t/n}(b)}
\le C_{n,t_0,R}
\norm[L^2]{a-b}.
\]

Therefore $w=q_n-v$ satisfies
\[
\norm[L^2]{w(t)}
\le C_{n,t_0,R}
\left|
\int_{t_0}^t\norm[L^2]{w(\tau)}
\,\dd\tau
\right|.
\]
Gronwall's inequality, forward and
backward from $t_0$, yields $w=0$ on
$I$. Hence $q_n=v\in C(I;H^3)$.

Since $t_0>0$ was arbitrary, and using
a one-sided neighborhood at $t=T_*$,
we conclude that
\begin{equation}\label{eq:positive-H3-continuity}
q_n\in C((0,T_*];H^3(\R)).
\end{equation}
\noindent
\textbf{Step 2. Weak $H^3$ continuity at zero.}
We already know that
\[
q_n(t)\longrightarrow u_0
\quad\text{strongly in }L^2(\R)
\quad\text{as }t\downarrow0,
\]
and that
\[
\sup_{0\le t\le T_*}
\norm[H^3]{q_n(t)}\le 4Y_0.
\]

We claim that
\begin{equation}\label{eq:weak-H3-at-zero}
q_n(t)\rightharpoonup u_0
\quad\text{weakly in }H^3(\R)
\quad\text{as }t\downarrow0.
\end{equation}

Indeed, let $t_k\downarrow0$ be arbitrary.
The sequence $(q_n(t_k))$ is bounded
in the Hilbert space $H^3$. Thus every
subsequence has a further subsequence
converging weakly in $H^3$ to some
$v\in H^3$.

Since $H^3$ embeds continuously into
$L^2$, the same subsequence converges
weakly in $L^2$ to $v$. However,
$q_n(t_k)\to u_0$ strongly in $L^2$.
Hence $v=u_0$.

Every weak cluster point in $H^3$ is
therefore equal to $u_0$, which proves
\eqref{eq:weak-H3-at-zero}.

\medskip
\noindent
\textbf{Step 3. Convergence of the $H^3$ norms.}
For each fixed $t$, the convergence
$q_{n,\varepsilon}(t)\to q_n(t)$
in $L^2$, together with the uniform
$H^3$ bound, implies weak $H^3$
convergence along a subsequence.
Thus \eqref{eq:Ybound} and the weak lower semicontinuity of the $H^3$-norm yields
\[
\norm[H^3]{q_n(t)}
\le
\frac{Y_0}{
\left(1-\frac C2t\sqrt{M_0Y_0}\right)^2}.
\]
Consequently,
\[
\limsup_{t\downarrow0}
\norm[H^3]{q_n(t)}
\le Y_0.
\]

On the other hand,
\eqref{eq:weak-H3-at-zero} and weak
lower semicontinuity imply
\[
Y_0=\norm[H^3]{u_0}
\le
\liminf_{t\downarrow0}
\norm[H^3]{q_n(t)}.
\]
Combining the two inequalities, it follows
\begin{equation}\label{eq:H3-norm-convergence}
\lim_{t\downarrow0}
\norm[H^3]{q_n(t)}
=\norm[H^3]{u_0}.
\end{equation}
\noindent
\textbf{Step 4. Strong continuity at zero.}
Since $H^3(\R)$ is a Hilbert space,
\begin{align*}
\norm[H^3]{q_n(t)-u_0}^2
={}&\norm[H^3]{q_n(t)}^2
+\norm[H^3]{u_0}^2-2
\ip{q_n(t)}{u_0}_{H^3}.
\end{align*}
By \eqref{eq:weak-H3-at-zero},
$
\lim_{t\downarrow0}\ip{q_n(t)}{u_0}_{H^3}
=\norm[H^3]{u_0}^2,
$ 
while \eqref{eq:H3-norm-convergence}
yields
$
\lim_{t\downarrow0}\norm[H^3]{q_n(t)}^2
=\norm[H^3]{u_0}^2.
$ 
Hence
$
\lim_{t\downarrow0}
\norm[H^3]{q_n(t)-u_0}=0.
$
{Combining this with}
\eqref{eq:positive-H3-continuity}, yields 
$
q_n\in C([0,T_*];H^3(\R)),
$
i.e. \eqref{eq:solutionreg}.
}
\end{proof}

\section{Well-posedness of the KdV equation}

\subsection{Consistency of the nonlinear split form}
\begin{lemma}[Consistency of the nonlinear split form]\label{lem:consistency}
For every real $f\in H^3$,
\begin{equation}\label{eq:cons2}
 \norm[L^2]{\F_h(f)-6fDf}
 \le Ch^2\norm[H^3]{f}^2.
\end{equation}
\end{lemma}
\begin{proof}
It follows from
$
 \F_h(f)-6fDf=2f(\delta_h-D)f+2(\delta_h-D)(f^2),
$
the multiplier inequality
$\norm[L^2]{(\delta_h-D)g}\le Ch^2\norm[H^3]{g}$,
and the algebra property of $H^3$.
\end{proof}

\subsection{Mesh-uniform $H^2$ one-sided stability}
The crucial estimate is the following mesh-uniform one-sided $H^2$
inequality. {
We put $\mathcal A_h=\delta_h$ 
} and write
\[
 G_h(q):=2q\mathcal A_hq+2\mathcal A_h(q^2).
\]
{
The operator $\mathcal A_h$ is skew-adjoint, commutes with $D$,}
and satisfies
\begin{equation}\label{eq:H2-Aprops} \multilabel{eq:delta-Lp-derivative}
 \|\mathcal A_hf\|_{L^p}\le \|Df\|_{L^p}
 \quad(1\leq p\leq \infty),\qquad
 \|[a,\mathcal A_h]f\|_2\le
  \|Da\|_\infty\|f\|_2.
\end{equation}
For $h=0$, the identities below are justified on the 
$H^3$ inputs (or first for smooth functions and then by density).

\begin{lemma}[Mesh-uniform $H^2$ one-sided stability]
\label{lem:H2onesided}
For real $q,p\in H^3(\R)$, set $a=q+p$, $w=q-p$.
There is a constant $C>0$ independent of $h$ such that
\begin{equation}\label{eq:H2onesided}
 \left|\sum_{k=0}^2
 \langle D^k(G_h(q)-G_h(p)),D^kw\rangle_{L^2}\right|
 \le C \|a\|_{H^3}\|w\|_{H^2}^2.
\end{equation}
\end{lemma}
\begin{proof}
Polarization gives the exact identity
\begin{equation}\label{eq:H2polar}
 G_h(q)-G_h(p)
 =a\mathcal A_hw+w\mathcal A_ha+2\mathcal A_h(aw).
\end{equation}
Fix $k\in\{0,1,2\}$ and put $f=D^kw$. By Leibniz' rule,
\begin{align*}
 D^k(G_h(q)-G_h(p))
 ={}&a\mathcal A_hf+2\mathcal A_h(af)
       +(\mathcal A_ha)f\\
 &+\sum_{i=1}^{k}\binom{k}{i}
 \bigl[(D^ia)\mathcal A_hD^{k-i}w
       +(D^{k-i}w)\mathcal A_hD^ia
 +2\mathcal A_h((D^ia)D^{k-i}w)\bigr].
\end{align*}
Skew-adjointness cancels the dangerous highest-order term:
\[
 \langle a\mathcal A_hf+2\mathcal A_h(af),f\rangle
 =-\langle a\mathcal A_hf,f\rangle.
\]
The commutator bound in \eqref{eq:H2-Aprops} controls this by
${
\frac12}\|Da\|_\infty\|f\|_2^2$.
The term $(\mathcal A_ha)f$ is bounded by
$ \|Da\|_\infty\|f\|_2^2$.

For the remaining terms, write $j=k-i$, so $i\ge1$ and
$i+j\le2$. From \eqref{eq:H2-Aprops} and the product rule,
\begin{align*}
 \|(D^ia)\mathcal A_hD^jw\|_2
 &\le \|D^ia\|_\infty\|D^{j+1}w\|_2,\\
 \|(D^jw)\mathcal A_hD^ia\|_2
 &\le \min\bigl\{\|D^jw\|_2\|D^{i+1}a\|_\infty,\;
                      \|D^jw\|_\infty\|D^{i+1}a\|_2\bigr\},\\
 \|\mathcal A_h((D^ia)D^jw)\|_2
 &\le \bigl(\|D^{i+1}a\,D^jw\|_2
                      +\|D^ia\,D^{j+1}w\|_2\bigr).
\end{align*}
In one dimension $H^3\hookrightarrow W^{2,\infty}$ and
$H^2\hookrightarrow W^{1,\infty}$. 
If $i=2$, $j=0$, use the second alternative in the middle line and
$\|D^3a\,w\|_2\le\|D^3a\|_2\|w\|_\infty$ in the last line;
if $i=1$, use the first alternative, so that only derivatives of $a$
of order at most two occur in the $L^\infty$ factors. Hence each remainder is bounded by
$C \|a\|_{H^3}\|w\|_{H^2}$ in $L^2$.
Pairing with $f=D^kw$, and summing over $k\le2$, proves
\eqref{eq:H2onesided}.
\end{proof}

\subsection{Cauchy property and construction of the KdV flow}
\begin{theorem}
[Cauchy property and construction of the KdV flow]
\label{thm:order2}
Let $u_0\in H^3(\R;\R)$ and let
$\{q_n\}_{n\ge1}$ be the solutions of
\eqref{eq:scheme} constructed in
Theorem~\ref{thm:wellposed}.

There exists $T>0$, depending only on
$\|u_0\|_{H^3}$, such that the following
assertions hold.

\begin{enumerate}[label=\textup{(\roman*)},
leftmargin=*]

\item The sequence $(q_n)_{n\ge1}$ is
Cauchy in $C([0,T];H^2(\R))$.
Consequently, there exists
\[
u\in C([0,T];H^2(\R))
   \cap L^\infty([0,T];H^3(\R))
\]
such that
\begin{equation}\label{eq:cauchy-limit}
q_n\longrightarrow u
\quad\text{in }C([0,T];H^2(\R)).
\end{equation}

\item 
{
 The limit $u$ satisfies
\[
u\in C([0,T];H^2(\R))
\cap C^1([0,T];H^{-1}(\R)),
\]
and, for $t\in[0,T]$,
\begin{equation}\label{eq:constructed-kdv}
\begin{cases}
u_t=-D^3u+6uDu
\quad\text{in }H^{-1}(\R),\\
u(0)=u_0.
\end{cases}
\end{equation}

Equivalently,
\[
u(t)=u_0+
\int_0^t
\bigl(-D^3u(\tau)+6u(\tau)Du(\tau)\bigr)
\,\dd\tau, \qquad t\in[0,T],
\]
where the integral is understood in the
Bochner sense in $H^{-1}(\R)$.
}

Moreover
\begin{equation}\label{eq:limit-mass}
\|u(t)\|_{L^2}=\|u_0\|_{L^2},
\qquad 0\le t\le T.
\end{equation}

\item If, in addition, $u_0\in H^7(\R)$,
then there exists $C_{T,u_0}>0$,
independent of $n,m$, such that
\begin{equation}\label{eq:two-mesh-cauchy}
\|q_n-q_m\|_{C_tH^2}
\le C_{T,u_0}|n^{-2}-m^{-2}|,
\qquad n,m\ge1.
\end{equation}
Consequently,
\begin{equation}\label{eq:rate2}
\|q_n-u\|_{C_tH^2}
\le C_{T,u_0}n^{-2}.
\end{equation}
\end{enumerate}
\end{theorem}

\begin{proof}
We divide the proof into six steps.

\medskip
\noindent
\textbf{Step 1. Uniform propagation of higher
Sobolev regularity.}

We first establish that, for every integer
$r\ge3$, the solutions constructed in
Theorem~\ref{thm:wellposed} satisfy
\begin{equation}\label{eq:high-reg-propagation}
\sup_{n\ge1}\sup_{0\le t\le T}
\|q_n(t)\|_{H^r}
\le C_{r,T}\|u_0\|_{H^r},
\end{equation}
provided $u_0\in H^r$.
The time $T$ depends only on the
$H^3$ bound of the initial datum.

To prove this assertion, we establish
the mesh-uniform estimate
\begin{equation}\label{eq:high-reg-energy}
\left|
\langle \F_h(f),f\rangle_{H^r}
\right|
\le
C_r\|Df\|_{L^\infty}
\|f\|_{H^r}^2,
\end{equation}
initially for real $f\in\mathcal S(\R)$,
where $C_r$ is independent of $h$.

Recall that
\[
\F_h(f)=2f\delta_hf+2\delta_h(f^2).
\]
Fix $1\le k\le r$ and put $f_k=D^kf$.
Since $D$ commutes with $\delta_h$,
Leibniz' rule gives
\begin{align*}
D^k\F_h(f)
={}&2f\delta_hf_k+4\delta_h(ff_k)+2\sum_{j=1}^{k}
\binom{k}{j}
(D^jf)\delta_hD^{k-j}f+2\delta_h\left(
\sum_{j=1}^{k-1}
\binom{k}{j}
(D^jf)(D^{k-j}f)
\right).
\end{align*}

The first two terms contain the
highest-order derivatives. By the
skew-adjointness of $\delta_h$,
\begin{align*}
&\left\langle
2f\delta_hf_k+4\delta_h(ff_k),
f_k\right\rangle_{L^2}
=-2\langle f\delta_hf_k,f_k\rangle_{L^2}.
\end{align*}
Lemma~\ref{lem:commutator} therefore yields
\[
\left|
\left\langle
2f\delta_hf_k+4\delta_h(ff_k),
f_k\right\rangle_{L^2}
\right|
\le
\|Df\|_\infty\|f_k\|_2^2.
\]

We next estimate the remaining terms.

For clarity, we recall the one-dimensional
Gagliardo--Nirenberg product inequality:
if $v\in L^\infty\cap H^\ell$,
$\ell\ge1$, and $\alpha+\beta=\ell$,
then\footnote{Indeed, apply the Gagliardo--Nirenberg
inequalities
\[
\|D^jv\|_{L^{2\ell/j}}
\le C_\ell
\|v\|_\infty^{1-j/\ell}
\|D^\ell v\|_2^{j/\ell}
\]
for $0<j\le\ell$, with the usual
$L^\infty$ interpretation at $j=0$,
and then use Holder's inequality.}
\begin{equation}\label{eq:GN-product-cauchy}
\|D^\alpha v\,D^\beta v\|_2
\le
C_\ell\|v\|_\infty\|D^\ell v\|_2.
\end{equation}

Apply \eqref{eq:GN-product-cauchy}
with $v=Df$ and $\ell=k-1$.
{
For $1\le j\le k-1$, put $v=Df$, $\ell=k-1$, and choose
$p=2\ell/(j-1)$ and $p'=2\ell/(k-j)$ (with $p=\infty$ when $j=1$).
Then $1/p+1/p'=1/2$, and the $L^p$ centered-difference bound yields
\begin{align*}
\|(D^jf)\delta_hD^{k-j}f\|_2
&\le\|D^{j-1}v\|_{L^p}
\|\delta_hD^{k-j-1}v\|_{L^{p'}}\\
&\le\|D^{j-1}v\|_{L^p}
\|D^{k-j}v\|_{L^{p'}}
\le C_k\|Df\|_\infty\|D^kf\|_2,
\end{align*}
by the one-dimensional Gagliardo--Nirenberg inequalities.
}
The endpoint $j=k$ satisfies the same
estimate, since
\[
\|D^kf\,\delta_hf\|_2
\le
\|D^kf\|_2\|Df\|_\infty.
\]

Likewise, for $1\le j\le k-1$,
\begin{align*}
&\|\delta_h((D^jf)(D^{k-j}f))\|_2
\le
\|D((D^jf)(D^{k-j}f))\|_2
\le
C_k\|Df\|_\infty\|D^kf\|_2,
\end{align*}
where the last estimate follows by
expanding the derivative and applying
\eqref{eq:GN-product-cauchy} to each
product.

The case $k=1$ is immediate from
$\|\delta_hf\|_\infty\le\|Df\|_\infty$.

Consequently,
\[
\left|
\langle D^k\F_h(f),D^kf\rangle_2
\right|
\le
C_k\|Df\|_\infty\|f\|_{H^k}^2.
\]
For $k=0$, the contribution vanishes
 by the skew-symmetric splitting.
Summing over $0\le k\le r$ proves
\eqref{eq:high-reg-energy}.
Density gives the corresponding
quadratic-form estimate on $H^r$.

Apply this estimate first to the
positive-mesh regularized solutions
$q_{n,\varepsilon}$.
The interaction-picture energy identity
of Lemma~\ref{lem:mildenergy}, whose
proof applies verbatim to every integer
Sobolev order, gives
\[
\frac{\dd}{\dd t}
\|q_{n,\varepsilon}(t)\|_{H^r}
\le
C_r\|Dq_{n,\varepsilon}(t)\|_\infty
\|q_{n,\varepsilon}(t)\|_{H^r}.
\]

By the uniform $H^3$ estimate,
\[
\|Dq_{n,\varepsilon}(t)\|_\infty
\le
C\|q_{n,\varepsilon}(t)\|_{H^3}
\le CR
\]
on $[0,T]$, where
\[
R=4\|u_0\|_{H^3}.
\]
Therefore Gronwall's inequality gives
\begin{equation}\label{eq:high-reg-regbound}
\|q_{n,\varepsilon}(t)\|_{H^r}
\le
\|u_0\|_{H^r}\exp(C_rRt)\le \|u_0\|_{H^r}\exp(C_rRT).
\end{equation}

For $u_0\in H^r$, the regularized
solutions are initially constructed in
$H^r$ by the same positive-mesh
contraction argument. The bound above
permits continuation to the common
interval $[0,T]$.

Passing to the limit $\varepsilon\downarrow0$
using the already established strong
$C_tL^2$ convergence and weak lower
semicontinuity in $H^r$ proves
\eqref{eq:high-reg-propagation}.
The positive-time local-solvability
argument, followed by weak-to-strong
continuity at $t=0$ using
\eqref{eq:high-reg-regbound}, also yields
$q_n\in C([0,T];H^r)$.

In particular, if $u_0\in H^7$, then
\begin{equation}\label{eq:M7-cauchy}
M_7:=
\sup_{n\ge1}\sup_{0\le t\le T}
\|q_n(t)\|_{H^7}<\infty.
\end{equation}

\medskip
\noindent
\textbf{Step 2. Quantitative comparison of the
linear generators at two meshes.}

For $h\ge0$, let
\[
\mathcal E_h=Q_h+h\partial_hQ_h,
\qquad \mathcal E_0=L.
\]
Then
\[
E_n(t)=\mathcal E_{t/n}.
\]

We prove the multiplier estimate
\begin{equation}\label{eq:E-two-mesh}
\|(\mathcal E_a-\mathcal E_b)f\|_{H^s}
\le
C|a^2-b^2|\|f\|_{H^{s+5}},
\qquad a,b\ge0.
\end{equation}

Recall from Section~1 that
\[
g(z)=\sin^3z+\frac12\sin^5z,
\]
and define
\[
H(z)=-2g(z)+zg'(z).
\]
The Fourier symbol of $\mathcal E_h$
is $ih^{-3}H(h\xi)$.

Differentiating with respect to $h$,
we obtain
\[
\partial_h
\left(h^{-3}H(h\xi)\right)
=
h^{-4}J(h\xi),
\]
where
\[
J(z)=-3H(z)+zH'(z)
=6g(z)-4zg'(z)+z^2g''(z).
\]

Since
\[
g(z)=z^3-\frac{37}{120}z^7+O(z^9),
\]
the terms of degree three cancel,
and
\[
J(z)=O(z^7)\qquad(z\to0).
\]
Furthermore, $g$, $g'$, and $g''$
are bounded on $\R$, so
\[
|J(z)|\le C(1+|z|^2)
\qquad (|z|\ge1).
\]
Combining the small- and
large-frequency estimates gives
\[
|J(z)|\le C|z|^5,
\qquad z\in\R.
\]

Consequently,
\[
\left|
\partial_h\left(h^{-3}H(h\xi)\right)
\right|
\le Ch|\xi|^5.
\]
Integration between $a$ and $b$
therefore yields
\[
|\mathcal E_a(\xi)-\mathcal E_b(\xi)|
\le C|a^2-b^2||\xi|^5.
\]
Then, Plancherel's theorem yields \eqref{eq:E-two-mesh}.
In particular,
\begin{equation}\label{eq:E-nm}
\|(E_n(t)-E_m(t))f\|_{H^2}
\le
Ct^2|n^{-2}-m^{-2}|
\|f\|_{H^7}.
\end{equation}

\medskip
\noindent
\textbf{Step 3. Quantitative comparison of
the nonlinearities at two meshes.}

By Lemma~\ref{lem:meshvariation},
\[
\|(\delta_a-\delta_b)f\|_{H^s}
\le
C|a^2-b^2|\|f\|_{H^{s+3}}
\]
for every $s\in\R$ and $a,b\ge0$.

Using
\[
(\F_a-\F_b)(f)
=
2f(\delta_a-\delta_b)f
+2(\delta_a-\delta_b)(f^2),
\]
and the Sobolev algebra properties
of $H^2$ and $H^5$, we obtain
\begin{equation}\label{eq:F-two-mesh-H2}
\|(\F_a-\F_b)(f)\|_{H^2}
\le
C|a^2-b^2|\|f\|_{H^5}^2.
\end{equation}

Indeed,
\begin{align*}
&\|f(\delta_a-\delta_b)f\|_{H^2}
\le
C\|f\|_{H^2}
\|(\delta_a-\delta_b)f\|_{H^2}
\le
C|a^2-b^2|\|f\|_{H^5}^2,
\end{align*}
and
\begin{align*}
\|(\delta_a-\delta_b)(f^2)\|_{H^2}
&\le
C|a^2-b^2|\|f^2\|_{H^5}\le
C|a^2-b^2|\|f\|_{H^5}^2.
\end{align*}

Consequently,
\begin{align}
&\|\F_{t/n}(f)-\F_{t/m}(f)\|_{H^2}
\le
Ct^2|n^{-2}-m^{-2}|\|f\|_{H^5}^2.
\label{eq:F-nm}
\end{align}

\medskip
\noindent
\textbf{Step 4. The Cauchy estimate for
smooth initial data.}

Assume for the moment that
$u_0\in H^7(\R)$.

Let
\[
w_{n,m}=q_n-q_m.
\]
Subtracting the equations satisfied
by $q_n$ and $q_m$, we find
\begin{align*}
\partial_tw_{n,m}
={}&E_n(t)w_{n,m}
+\F_{t/n}(q_n)-\F_{t/n}(q_m)
+(E_n(t)-E_m(t))q_m
+\F_{t/n}(q_m)-\F_{t/m}(q_m).
\end{align*}

Thus
\begin{equation}\label{eq:cauchy-difference}
\partial_tw_{n,m}
=
E_n(t)w_{n,m}
+\F_{t/n}(q_n)-\F_{t/n}(q_m)
+\mathcal R_{n,m},
\end{equation}
where
\begin{align}
\mathcal R_{n,m}(t)
={}&(E_n(t)-E_m(t))q_m(t)
+\F_{t/n}(q_m(t))
-\F_{t/m}(q_m(t)).
\label{eq:cauchy-remainder}
\end{align}

By \eqref{eq:M7-cauchy},
\eqref{eq:E-nm}, and \eqref{eq:F-nm},
we obtain
\begin{equation}\label{eq:cauchy-forcing}
\|\mathcal R_{n,m}(t)\|_{H^2}
\le
C(M_7+M_7^2)
t^2|n^{-2}-m^{-2}|.
\end{equation}

{Since $q_n,q_m\in C_tH^3\cap C_t^1L^2$ satisfy their differential equations in $L^2$, equation \eqref{eq:cauchy-difference} holds in $L^2$. Conjugating by $U_n(0,t)$ and applying Proposition~\ref{prop:differential-mild} gives the following variation-of-constants formula (the forcing is continuous in $H^2$ for $H^7$ data).}
Since $w_{n,m}(0)=0$, the corresponding
mild equation is
\begin{align*}
w_{n,m}(t)
=\int_0^t U_n(t,\tau)
\Bigl[
&\F_{\tau/n}(q_n(\tau))
-\F_{\tau/n}(q_m(\tau))+\mathcal R_{n,m}(\tau)
\Bigr]\,\dd\tau.
\end{align*}

The forcing belongs to $L^1([0,T];H^2)$
by the $H^7$ regularity established
in Step~1. Hence the interaction-picture
energy identity of
Lemma~\ref{lem:mildenergy} applies
with $k=2$.

Using the mesh-uniform one-sided $H^2$
estimate of Lemma~\ref{lem:H2onesided}
and the uniform $H^3$ bound, we get,
for almost every $t\in[0,T]$,
\begin{align*}
\frac12\frac{\dd}{\dd t}
\|w_{n,m}(t)\|_{H^2}^2
\le{}
CR\|w_{n,m}(t)\|_{H^2}^2
+\|\mathcal R_{n,m}(t)\|_{H^2}
\|w_{n,m}(t)\|_{H^2}.
\end{align*}

Applying the standard regularization
of the norm, or equivalently the
integral form of Gronwall's inequality,
we obtain
\[
\|w_{n,m}(t)\|_{H^2}
\le
\int_0^t
e^{CR(t-\tau)}
\|\mathcal R_{n,m}(\tau)\|_{H^2}
\,\dd\tau.
\]

By \eqref{eq:cauchy-forcing},
\begin{align*}
\|w_{n,m}(t)\|_{H^2}
&\le
C(M_7+M_7^2)
|n^{-2}-m^{-2}|
\int_0^t e^{CR(t-\tau)}\tau^2
\,\dd\tau.
\end{align*}

Since $0\le t\le T$, we conclude that
\begin{equation}\label{eq:smooth-cauchy}
{
\|q_n-q_m\|_{C_tH^2}
\le C_{T,u_0}|n^{-2}-m^{-2}|.
}
\end{equation}

In particular, $(q_n)$ is Cauchy
in $C([0,T];H^2)$ whenever
$u_0\in H^7$.

\medskip
\noindent
\textbf{Step 5. Extension to arbitrary
$H^3$ initial data by density.}

Let now $0\not=u_0\in H^3$.
%
{
Choose an even, nonnegative $\rho\in C_c^\infty(\R)$ with
$\int_\R\rho=1$ and put
$u_0^{(j)}=\rho_{1/j}*u_0$.
Then $u_0^{(j)}\in H^r$ for every finite $r$,
$u_0^{(j)}\to u_0$ in $H^3$, and
\[
\|u_0^{(j)}\|_2\le\|u_0\|_2,\qquad
\|u_0^{(j)}\|_{H^3}\le\|u_0\|_{H^3}.
\]
}
Denote by $q_n^{(j)}$ the solution
of \eqref{eq:scheme} with initial
datum $u_0^{(j)}$. Then, the lifespan of $q_n^{(j)}$ $T_*(u_0^{(j)})\overset{\eqref{eq:Ybound}}\ge T_*(u_0)=: T$.

For each fixed $j$, Step~4 gives
\begin{equation}\label{eq:smooth-cauchy-j}
\|q_n^{(j)}-q_m^{(j)}\|_{C_tH^2}
\le
C_j|n^{-2}-m^{-2}|,\qquad t\in [0,T].
\end{equation}

We next use continuous dependence
on the initial data, uniformly in $n$.

Because $q_n$ and $q_n^{(j)}$ solve
the same equation with the same mesh,
Lemma~\ref{lem:H2onesided} and the
interaction-picture energy identity
give
\begin{equation}\label{eq:initial-stability-H2}
\|q_n-q_n^{(j)}\|_{C_tH^2}
\le
e^{CR_0T}
\|u_0-u_0^{(j)}\|_{H^2},
\end{equation}
where $C$ is independent of $n,j$.

By the triangle inequality,
\begin{align*}
\|q_n-q_m\|_{C_tH^2}
\le{}&
\|q_n-q_n^{(j)}\|_{C_tH^2}
+\|q_n^{(j)}-q_m^{(j)}\|_{C_tH^2}
+\|q_m^{(j)}-q_m\|_{C_tH^2}.
\end{align*}
Therefore
\begin{align}
\|q_n-q_m\|_{C_tH^2}
\le{}&
2e^{CR_0T}
\|u_0-u_0^{(j)}\|_{H^2}
+C_j|n^{-2}-m^{-2}|.
\label{eq:density-cauchy}
\end{align}

Fix $j$ and let $n,m\to\infty$.
We obtain
\[
\limsup_{n,m\to\infty}
\|q_n-q_m\|_{C_tH^2}
\le
2e^{CR_0T}
\|u_0-u_0^{(j)}\|_{H^2}.
\]

Finally, letting $j\to\infty$ gives
\[
\lim_{n,m\to\infty}
\|q_n-q_m\|_{C_tH^2}=0.
\]

Thus $(q_n)$ is Cauchy in
$C([0,T];H^2)$ for every
$u_0\in H^3$.

By completeness, there exists
\[
u\in C([0,T];H^2)
\]
such that
\begin{equation}\label{eq:q-to-u-H2}
q_n\longrightarrow u
\quad\text{in }C([0,T];H^2).
\end{equation}

Moreover, the uniform $H^3$ bound
and weak lower semicontinuity imply
\[
u(t)\in H^3,\qquad
\sup_{0\le t\le T}\|u(t)\|_{H^3}
\le R.
\]
Thus $u\in L^\infty_tH^3\cap C_tH^2$.

\medskip
\noindent
\textbf{Step 6. Identification of the limit
as a KdV solution.}

We now show that the limit $u$
satisfies \eqref{eq:constructed-kdv}.

For every $n$, the differential
formulation holds in $L^2$:
\[
q_{n,t}
=E_n(t)q_n+\F_{t/n}(q_n).
\]
Consequently,
\begin{equation}\label{eq:q-integrated}
q_n(t)=u_0+
\int_0^t
\left[
E_n(\tau)q_n(\tau)
+\F_{\tau/n}(q_n(\tau))
\right]\,\dd\tau.
\end{equation}

We can pass to the limit in
$C([0,T];H^{-1})$.Indeed:

\smallskip
\noindent
\emph{$\bullet$ Convergence of the linear term.}
The multiplier estimate \eqref{eq:Ebound}
with $r=-1$ gives
\[
\|E_n(t)f\|_{H^{-1}}
\le C\|f\|_{H^2}
\]
uniformly in $n,t$.
Therefore,
\[
\|E_n(t)(q_n(t)-u(t))\|_{H^{-1}}
\le C\|q_n(t)-u(t)\|_{H^2},
\]
which converges to zero uniformly in $t$.

Moreover, $E_n(t)f\to Lf$ in
$H^{-1}$, uniformly for $t\in[0,T]$,
for every fixed $f\in H^2$.
This follows by proving the assertion
first for smooth $f$ using
\eqref{eq:Econs}, and then applying
density together with the uniform
$H^2\to H^{-1}$ bounds.
Since
$
\{u(t):0\le t\le T\}
$
is compact in $H^2$, the convergence
is uniform when $f=u(t)$.
Consequently,
\begin{equation}\label{eq:linear-limit-cauchy}
E_n(t)q_n(t)\longrightarrow Lu(t)
\quad\text{in }C([0,T];H^{-1}).
\end{equation}

\smallskip
\noindent
\emph{$\bullet$ Convergence of the nonlinear term.}
The centered difference satisfies
\[
\|\delta_hf\|_{H^1}
\le\|f\|_{H^2},
\]
uniformly in $h$.

Using the algebra properties of
$H^1$ and $H^2$, we obtain
\[
\|\F_h(f)-\F_h(g)\|_{H^1}
\le
C(\|f\|_{H^2}+\|g\|_{H^2})
\|f-g\|_{H^2},
\]
with $C$ independent of $h$.

Hence \eqref{eq:q-to-u-H2} implies
\[
\F_{t/n}(q_n(t))-\F_{t/n}(u(t))
\longrightarrow0
\quad\text{in }C_tH^1.
\]

For every fixed $f\in H^2$,
\[
\delta_hf\longrightarrow Df
\quad\text{in }H^1
\qquad(h\to0).
\]
By the same compactness argument
used for the linear term, this
convergence is uniform for
$f=u(t)$, $0\le t\le T$.

It follows that
\[
\F_{t/n}(u(t))
\longrightarrow
2u(t)Du(t)+2D(u(t)^2)
=6u(t)Du(t)
\]
in $C([0,T];H^1)$.

Therefore
\begin{equation}\label{eq:nonlinear-limit-cauchy}
\F_{t/n}(q_n(t))
\longrightarrow6u(t)Du(t)
\quad\text{in }C([0,T];H^1).
\end{equation}

Passing to the limit in
\eqref{eq:q-integrated}, using
\eqref{eq:linear-limit-cauchy} and
\eqref{eq:nonlinear-limit-cauchy},
we obtain
\[
u(t)=u_0+
\int_0^t
\bigl[-D^3u(\tau)+6u(\tau)Du(\tau)\bigr]
\,\dd\tau
\]
in $H^{-1}$.

The integrand is continuous in
$H^{-1}$, and thus
\[
u\in C^1([0,T];H^{-1}),
\]
with
\[
u_t=-D^3u+6uDu.
\]
This proves \eqref{eq:constructed-kdv}.

Finally, every $q_n$ conserves the
$L^2$ norm:
\[
\|q_n(t)\|_2=\|u_0\|_2.
\]
The strong $C_tH^2$ convergence
implies strong $C_tL^2$ convergence.
Therefore
\[
\|u(t)\|_2=\|u_0\|_2,
\]
which proves \eqref{eq:limit-mass}.

Uniqueness is proven in Lemma~\ref{lem:kdv-H2-uniqueness} below.

Finally, when $u_0\in H^7$,
letting $m\to\infty$ in
\eqref{eq:smooth-cauchy} and using
\eqref{eq:q-to-u-H2} gives
\[
\|q_n-u\|_{C_tH^2}
\le C_{T,u_0}n^{-2}.
\]

This completes the proof.
\end{proof}

\subsection{Lipschitz dependence of the KdV flow on the data}
{
\begin{corollary}[Lipschitz Regularity of the KdV flow]
\label{cor:kdv-H2-continuous-dependence}
{
For every $\rho>0$ there is $T(\rho)>0$ such that, for all real $f,g\in H^3(\R)$ with
$\|f\|_{H^3},\|g\|_{H^3}\le\rho$, then
\begin{equation}\label{eq:constructed-flow-lipschitz}
\|S(\cdot)f-S(\cdot)g\|_{C([0,T(\rho)];H^2)}
\le e^{C\rho T(\rho)}\|f-g\|_{H^2}.
\end{equation}
{
{Moreover,} 
there exists $C_\rho>0$ such that
\begin{equation}\label{eq:local-joint-lipschitz}
 \|S(t)f-S(s)g\|_{L^2}
 \le C_\rho(\|f-g\|_{H^2}+|t-s|),
 \quad t,s\in[0,T(\rho)].
\end{equation}
}
}
\end{corollary}
\begin{proof}
{
Fix $\rho>0$ and choose $T=T(\rho)>0$ so that
Theorem~\ref{thm:wellposed} supplies solutions $q_n^f,q_n^g$
for any $\|f\|_{H^3},\|g\|_{H^3}\le\rho$, with
$\sup_{n,t}(\|q_n^f(t)\|_{H^3}+\|q_n^g(t)\|_{H^3})\le C\rho$.
For a common $n$ let $w_n=q_n^f-q_n^g$.
The two functions have the same linear propagator $U_n$ and the same
nonlinear mesh $t/n$, so the interaction-picture identity and
Lemma~\ref{lem:H2onesided} give
$\frac12\frac{\dd}{\dd t}\|w_n\|_{H^2}^2
 \le C\rho\|w_n\|_{H^2}^2$ (a.e. in time).
The linear contribution vanishes by skew-adjointness.
Gronwall's inequality therefore yields
$\|q_n^f-q_n^g\|_{C_tH^2}
 \le e^{C\rho T}\|f-g\|_{H^2}$, with constant independent of $n$.
Theorem~\ref{thm:order2} supplies strong $C_tH^2$ convergence
for each initial datum; passing to the limit in the preceding
inequality proves \eqref{eq:constructed-flow-lipschitz}.
The factor $e^{C\rho T}$ is uniform on the $H^3$ ball, which is
precisely the claimed local Lipschitz dependence in $H^2$.
{By Proposition~\ref{prop:restart-H3},
$S(\cdot)g\in C_tH^3\cap C_t^1L^2$ and
$\partial_tS(t)g=-D^3S(t)g+6S(t)gD S(t)g$. The uniform $H^3$
bound on the ball implies $\|\partial_tS(t)g\|_2
\le C(\rho+\rho^2)$, independently of $g$ in that ball.
Thus $\|S(t)g-S(s)g\|_2\le C(\rho+\rho^2)|t-s|$.
Combining this with \eqref{eq:constructed-flow-lipschitz} and
$\|S(t)f-S(t)g\|_2\le\|S(t)f-S(t)g\|_{H^2}$
gives \eqref{eq:local-joint-lipschitz}.
}}
\end{proof}
}

\subsection{Uniqueness for KdV in $C_tH^2$}

\begin{lemma}[Uniqueness for KdV in $C_tH^2$]
\label{lem:kdv-H2-uniqueness}
Let
$
u,v\in C([0,T];H^2(\R))
\cap C^1([0,T];H^{-1}(\R))
$
be real-valued solutions of
$
u_t=-D^3u+6uDu
$
in $H^{-1}(\R)$.
Then, for every $t\in[0,T]$,
\begin{align}
\|u(t)-v(t)\|_{L^2}^2
\le{}&
\|u(0)-v(0)\|_{L^2}^2
\exp\left(
3\int_0^t
\|D(u(s)+v(s))\|_{L^\infty}\,ds
\right).
\label{eq:kdv-H2-stability}
\end{align}
In particular, solutions in this class
with the same initial datum coincide.
\end{lemma}

\begin{proof}
Set
$
w=u-v,\qquad a=u+v.
$
Since $6uDu=3D(u^2)$, subtraction gives
\[
w_t=-D^3w+3D(aw)
\quad\text{in }H^{-1}(\R).
\]

Since \(w\in C([0,T];H^2)\) and \(w_t\in C([0,T];H^{-1})\), the \(L^2\) energy identity, justified by spatial Fourier regularization, yields for almost every $t$,
\[
\frac12\frac{\dd}{\dd t}
\|w(t)\|_2^2
=
\langle w_t(t),w(t)\rangle_{H^{-1},H^1}.
\]

Since $D^3$ is skew-adjoint,
\[
\langle -D^3w,w\rangle_{H^{-1},H^1}
=
\int_{\R}(D^2w)(Dw)\,\dd x=0.
\]

For the nonlinear term,
\begin{align*}
3\langle D(aw),w\rangle
&=-3\int_{\R}awDw\,\dd x=\frac32\int_{\R}(Da)w^2\,\dd x.
\end{align*}
Consequently,
\[
\frac{\dd}{\dd t}\|w(t)\|_2^2
=
3\int_{\R}(Da)w^2\,\dd x
\]
for almost every $t\in[0,T]$.

By the Sobolev embedding
$H^2(\R)\hookrightarrow W^{1,\infty}(\R)$,
\[
\frac{\dd}{\dd t}\|w(t)\|_2^2
\le
3\|Da(t)\|_\infty\|w(t)\|_2^2.
\]

Gronwall's inequality yields
\eqref{eq:kdv-H2-stability}.
If $u(0)=v(0)$, the right-hand side
vanishes for every $t$, and therefore
$u=v$ on $[0,T]$.
\end{proof}

\subsection{Strong $H^3$ continuity by restarting the KdV flow}

\begin{proposition}
[Strong $H^3$ continuity]
\label{prop:restart-H3}
Let $u_0\in H^3(\R;\R)$, and let $u$
be the solution constructed in
Theorem~\ref{thm:order2}. Then
\[
u\in C([0,T];H^3(\R))
\cap C^1([0,T];L^2(\R)).
\]
\end{proposition}

\begin{proof}
We exploit the autonomy and reversibility
of KdV, together with 
{
uniqueness in $C_tH^2\cap C_t^1H^{-1}$, as proved in
Lemma~\ref{lem:kdv-H2-uniqueness}.}

\medskip
\noindent
\textbf{Step 1. Strong continuity at the
initial time for arbitrary data.}

Let $f\in H^3(\R;\R)$ be arbitrary.
Denote by $S(t)f$ the KdV solution
constructed by Theorem~\ref{thm:order2}
with initial datum $f$, and by
$q_n^f$ its Chernoff approximations.

The uniform energy estimate
\eqref{eq:Ybound} yields
\[
\|q_n^f(t)\|_{H^3}
\le
\frac{\|f\|_{H^3}}
{\bigl(1-c t
\sqrt{\|f\|_{L^2}\|f\|_{H^3}}\bigr)^2},
\]
on their common existence interval,
where $c>0$ is independent of $n$.

Since $q_n^f\to S(\cdot)f$ in $C_tH^2$,
weak lower semicontinuity in $H^3$
gives
\[
\|S(t)f\|_{H^3}
\le
\frac{\|f\|_{H^3}}
{\bigl(1-c t
\sqrt{\|f\|_{L^2}\|f\|_{H^3}}\bigr)^2}.
\]
Consequently,
\[
\limsup_{t\downarrow0}
\|S(t)f\|_{H^3}
\le \|f\|_{H^3}.
\]

On the other hand, strong $H^2$
continuity and the uniform $H^3$
bound imply
\[
S(t)f\rightharpoonup f
\quad\text{weakly in }H^3
\quad(t\downarrow0).
\]
Weak lower semicontinuity therefore gives
\[
\|f\|_{H^3}
\le
\liminf_{t\downarrow0}
\|S(t)f\|_{H^3}.
\]

Thus
\[
\lim_{t\downarrow0}
\|S(t)f\|_{H^3}
=\|f\|_{H^3}.
\]
Since $H^3$ is a Hilbert space,
weak convergence together with
convergence of norms implies
\begin{equation}\label{eq:restart-initial-H3}
\lim_{t\downarrow0}
\|S(t)f-f\|_{H^3}=0.
\end{equation}

This conclusion holds for every
$f\in H^3(\R;\R)$.

\medskip
\noindent
\textbf{Step 2. Right continuity at
arbitrary positive times.}

Fix $t_0\in[0,T)$ and put
$f=u(t_0)\in H^3$.

Since KdV is autonomous, the function
$v(t)=u(t_0+t)$ solves KdV with
initial datum $f$, on the interval
where it is defined.

By uniqueness of KdV solutions in
$C_tH^2\cap L^\infty_tH^3$, the
translated solution coincides locally
with the solution constructed from $f$:
$
u(t_0+t)=S(t)u(t_0).
$

Applying \eqref{eq:restart-initial-H3}
with $f=u(t_0)$ yields
\[
\lim_{t\downarrow0}
\|u(t_0+t)-u(t_0)\|_{H^3}=0.
\]
Hence $u$ is strongly right-continuous
in $H^3$ at every $t_0<T$.

\medskip
\noindent
\textbf{Step 3. Left continuity
by time-space reversal.}

Fix $t_0\in(0,T]$ and define
\[
v(s,x)=u(t_0-s,-x),
\qquad 0\le s\le t_0.
\]

The KdV equation is invariant under
the simultaneous reversal
$(t,x)\mapsto(-t,-x)$.
Indeed,
\[
v_s=-u_t,\qquad
v_{xxx}=-u_{xxx},\qquad
vv_x=-uu_x,
\]
where all quantities on the right
are evaluated at $(t_0-s,-x)$.
Therefore
\[
v_s=-v_{xxx}+6vv_x.
\]

In particular, $v$ solves KdV
with initial datum
\[
v(0,x)=u(t_0,-x)\in H^3.
\]

By uniqueness, $v$ coincides locally
with the solution obtained from its
initial datum through
Theorem~\ref{thm:order2}.
Applying \eqref{eq:restart-initial-H3}
to this initial datum gives
\[
\lim_{s\downarrow0}
\|v(s)-v(0)\|_{H^3}=0.
\]

Since spatial reflection preserves
the $H^3$ norm,
\[
\lim_{s\downarrow0}
\|u(t_0-s)-u(t_0)\|_{H^3}=0.
\]
Hence $u$ is strongly left-continuous
in $H^3$ at every $t_0>0$.

Combining the two one-sided
continuity statements, we conclude
\[
u\in C([0,T];H^3(\R)).
\]

\medskip
\noindent
\textbf{Step 4. Upgrading the
time regularity.}

Since $u\in C_tH^3$, we have
\[
-D^3u+6uDu\in C([0,T];L^2).
\]

The $H^{-1}$-valued integral equation
from Theorem~\ref{thm:order2}
therefore holds in $L^2$, by
injectivity of the embedding
$L^2\hookrightarrow H^{-1}$.

The fundamental theorem of calculus
for Bochner integrals yields
\[
u\in C^1([0,T];L^2),
\qquad
u_t=-D^3u+6uDu
\quad\text{in }L^2.
\]
\end{proof}

\subsection{Rate of convergence}

\begin{theorem}
[Rate of convergence]
\label{thm:kdv-sobolev-quantitative}
Let $3\le s\le7$ and let
$f\in H^s(\mathbb R;\mathbb R)$.
Let $q_n^f$ be the centered-difference approximants
of \eqref{eq:scheme}, and let $u=S(\cdot)f$
be the KdV flow.

For every $0<T\le T_*(f)$, there exists
$C_{T,f,s}>0$, independent of $n\ge1$, such that
\begin{equation}\label{eq:kdv-quant-H2}
{
\|q_n^f-u\|_{C([0,T];H^2)}
\le C_{T,f,s}n^{-\frac{2(s-2)}5}.
}
\end{equation}

Moreover, for every $f\in H^3(\mathbb R;\mathbb R)$,
\begin{equation}\label{eq:kdv-quant-L2}
{
\|q_n^f-u\|_{C([0,T];L^2)}
\le C_{T,f}n^{-6/5}.
}
\end{equation}
\end{theorem}

\begin{proof}
We sharpen the smooth-data two-mesh estimate
so that its constant depends linearly on the
high Sobolev norm. We then use Fourier
regularization of the initial datum.

All constants below are uniform in the
approximation parameters $n,m$, on a common
interval $[0,T]$. They may depend on $T$
and on the $H^3$-energy of the data.

\medskip
\noindent
\textbf{Step 1. A tame nonlinear mesh estimate.}

Recall that
\[
\mathcal F_h(v)
=2v\delta_hv+2\delta_h(v^2).
\]
By Lemma~\ref{lem:meshvariation},
\[
\|(\delta_a-\delta_b)w\|_{H^k}
\le C|a^2-b^2|\|w\|_{H^{k+3}},
\qquad k=0,2.
\]

For $v\in H^5$, the algebra property of $H^2$
and the tame product inequality
\[
\|v^2\|_{H^5}
\le C\|v\|_{H^1}\|v\|_{H^5}
\]
yields
\begin{align}
\|(\mathcal F_a-\mathcal F_b)(v)\|_{H^2}
&\le
2\|v(\delta_a-\delta_b)v\|_{H^2}
+2\|(\delta_a-\delta_b)(v^2)\|_{H^2}
\notag\\
&\quad\le
C|a^2-b^2|
\left(
\|v\|_{H^2}\|v\|_{H^5}
+\|v^2\|_{H^5}
\right)
\notag\\
&\quad\le
C|a^2-b^2|
\|v\|_{H^3}\|v\|_{H^5}.
\label{eq:kdv-tame-mesh-H2}
\end{align}

Similarly, the $H^3$ algebra property yields
\begin{equation}\label{eq:kdv-tame-mesh-L2}
\|(\mathcal F_a-\mathcal F_b)(v)\|_{L^2}
\le C|a^2-b^2|\|v\|_{H^3}^2.
\end{equation}

\medskip
\noindent
\textbf{Step 2. Tame smooth-data convergence estimates.}

Let $g\in H^7$, $\|g\|_{H^3}\le R$,
and $T\le T_*(g)$.

For $w_{n,m}=q_n^g-q_m^g$, the difference
equation \eqref{eq:cauchy-difference} has forcing
\[
\begin{aligned}
\mathcal R_{n,m}(t)
={}&(E_n(t)-E_m(t))q_m^g(t)+(\mathcal F_{t/n}-\mathcal F_{t/m})
(q_m^g(t)).
\end{aligned}
\]

The linear two-mesh estimate
\eqref{eq:E-nm} and
\eqref{eq:kdv-tame-mesh-H2} imply
\begin{align*}
\|\mathcal R_{n,m}(t)\|_{H^2}
&\le Ct^2|n^{-2}-m^{-2}|
\left(
\|q_m^g(t)\|_{H^7}
+\|q_m^g(t)\|_{H^3}
 \|q_m^g(t)\|_{H^5}
\right).
\end{align*}

By the uniform high-Sobolev propagation
estimate \eqref{eq:high-reg-propagation},
\[
\sup_{m,t}\|q_m^g(t)\|_{H^r}
\le C_{r,T,R}\|g\|_{H^r},
\qquad r=3,5,7.
\]
Consequently,
\[
\|\mathcal R_{n,m}(t)\|_{H^2}
\le C_{T,R}t^2|n^{-2}-m^{-2}|
\|g\|_{H^7}.
\]

Applying Lemma~\ref{lem:H2onesided}
and the interaction-picture energy argument
from Step~4 of Theorem~\ref{thm:order2},
we obtain
\[
\|q_n^g-q_m^g\|_{C_TH^2}
\le C_{T,R}|n^{-2}-m^{-2}|
\|g\|_{H^7}.
\]
Letting $m\to\infty$ yields
\begin{equation}\label{eq:kdv-smooth-tame-H2}
\|q_n^g-S(\cdot)g\|_{C_TH^2}
\le C_{T,R}n^{-2}\|g\|_{H^7}.
\end{equation}

We also derive a lower-regularity estimate
in $L^2$.

Let $g\in H^5$ with $\|g\|_{H^3}\le R$.
The linear estimate \eqref{eq:E-two-mesh},
with Sobolev index zero, gives
\[
\|(E_n(t)-E_m(t))v\|_2
\le Ct^2|n^{-2}-m^{-2}|
\|v\|_{H^5}.
\]
Combining it with
\eqref{eq:kdv-tame-mesh-L2} gives
\begin{align*}
\|\mathcal R_{n,m}(t)\|_2
&\le Ct^2|n^{-2}-m^{-2}|
\left(
\|q_m^g(t)\|_{H^5}
+\|q_m^g(t)\|_{H^3}^2
\right)\le C_{T,R}t^2|n^{-2}-m^{-2}|
\|g\|_{H^5}.
\end{align*}

The one-sided $L^2$ stability
estimate \eqref{eq:onesided}, together with
Lemma~\ref{lem:mildenergy} for $k=0$,
therefore yields
\[
\|q_n^g-q_m^g\|_{C_TL^2}
\le C_{T,R}|n^{-2}-m^{-2}|
\|g\|_{H^5}.
\]
Passing to the limit $m\to\infty$ gives
\begin{equation}\label{eq:kdv-smooth-tame-L2}
\|q_n^g-S(\cdot)g\|_{C_TL^2}
\le C_{T,R}n^{-2}\|g\|_{H^5}.
\end{equation}

\medskip
\noindent
\textbf{Step 3. Fourier regularization
and the $H^2$ convergence rate.}

Choose a real even function
$\chi\in C_c^\infty(\mathbb R)$ with
$0\le\chi\le1$ and $\chi=1$ on $[-1,1]$.

For $0<\delta\le1$, define
\[
\widehat{f_\delta}(\xi)
=\chi(\delta\xi)\widehat f(\xi).
\]

For $3\le s\le7$, we have
\begin{equation}\label{eq:kdv-cutoff-H2}
\begin{aligned}
\|f-f_\delta\|_{H^2}
&\le C\delta^{s-2}\|f\|_{H^s},\\
\|f_\delta\|_{H^7}
&\le C\delta^{s-7}\|f\|_{H^s}.
\end{aligned}
\end{equation}

Moreover,
\[
\|f_\delta\|_2\le\|f\|_2,
\qquad
\|f_\delta\|_{H^3}\le\|f\|_{H^3}.
\]
Hence the explicit lifespan formula gives
\[
T_*(f_\delta)\ge T_*(f).
\]

The uniform $H^2$ data-stability estimate
obtained from Lemma~\ref{lem:H2onesided}
and its limit for the KdV flow imply
\[
\begin{aligned}
&\|q_n^f-q_n^{f_\delta}\|_{C_TH^2}+
\|S(\cdot)f-S(\cdot)f_\delta\|_{C_TH^2}
\le C_{T,R}\|f-f_\delta\|_{H^2},
\end{aligned}
\]
where $R$ depends only on $\|f\|_{H^3}$.

Using the triangle inequality and
\eqref{eq:kdv-smooth-tame-H2}, we obtain
\begin{align*}
\|q_n^f-S(\cdot)f\|_{C_TH^2}
&\le C_{T,R}
\left(
\|f-f_\delta\|_{H^2}
+n^{-2}\|f_\delta\|_{H^7}
\right)\\
&\le C_{T,R}\|f\|_{H^s}
\left(
\delta^{s-2}
+n^{-2}\delta^{s-7}
\right).
\end{align*}

Choosing
$
{\delta=n^{-2/5},}
$
then
$
\delta^{s-2}
=n^{-2(s-2)/5}
=n^{-2}\delta^{s-7}.
$
which proves \eqref{eq:kdv-quant-H2}.

\medskip
\noindent
\textbf{Step 4. Improved $L^2$ rate
for $H^3$ initial data.}

For $f\in H^3$, the same frequency cutoff
satisfies
\[
\|f-f_\delta\|_2
\le C\delta^3\|f\|_{H^3},
\qquad
\|f_\delta\|_{H^5}
\le C\delta^{-2}\|f\|_{H^3}.
\]

The uniform $L^2$ data-stability estimate
of Corollary~\ref{cor:uniqueness}
also holds for the limiting KdV flow.

Using \eqref{eq:kdv-smooth-tame-L2},
we obtain
\begin{align*}
\|q_n^f-S(\cdot)f\|_{C_TL^2}
&\le C_{T,R}
\left(
\|f-f_\delta\|_2
+n^{-2}\|f_\delta\|_{H^5}
\right)\\
&\le C_{T,R}\|f\|_{H^3}
\left(
\delta^3+n^{-2}\delta^{-2}
\right).
\end{align*}

Again choosing $\delta=n^{-2/5}$ gives
$
\delta^3=n^{-6/5}
=n^{-2}\delta^{-2}.
$
Hence
\[
\|q_n^f-S(\cdot)f\|_{C_TL^2}
\le C_{T,f}n^{-6/5}.
\]
This completes the proof.
\end{proof}

\subsection{Proof of Theorem~\ref{prop:explicit-formula}}

\par\noindent{
\emph{Proof of Theorem~\ref{prop:explicit-formula}.}
For fixed $n,\varepsilon>0$ the positive-mesh contraction lemma yields the indicated
Picard limit on each successive interval (starting at $0$),
and the pieces concatenate to $p_{n,\varepsilon}$.
For fixed $n$, 
Proposition~\ref{prop:Cauchy} and Theorem~\ref{thm:order2} implies \eqref{eq:explicit-nonlinear-flow}.
{Indeed, for $f\in H^3(\R;\R)$, estimate \eqref{eq:reg-H2-quantitative} on $[0,T]$ reads
$\|p_{n,\varepsilon}-q_n\|_{C_TH^2}\le C_{T,f}\varepsilon n^{-2/3}$,
with a constant independent of $n,\varepsilon$. Hence
\[
 \|p_{n,\varepsilon}-S(\cdot)f\|_{C_TH^2}
 \le C_{T,f}\varepsilon n^{-2/3}
 +\|q_n-S(\cdot)f\|_{C_TH^2}\overset{n\to\infty}\longrightarrow0
\]
for every fixed $\varepsilon\in(0,T]$, by Theorem~\ref{thm:order2}.
This proves \eqref{eq:fixed-epsilon-H3}. 
}
Proposition~\ref{prop:restart-H3} proves \eqref{eq:KdV-Main}. 
%
$L^2$-norm conservation for $p_{n,\varepsilon}$ and $q_n$
follows from Proposition~\ref{prop:mass} and passage to the limit;
conservation for $S$ follows from strong $C_tL^2$ convergence.

{
\emph{Rates for $p_{n,\varepsilon}$.}
By \eqref{eq:reg-quantitative}, $\|p_{n,\varepsilon}-q_n\|_{C_TL^2}\le C_{T,f}\varepsilon^3n^{-2}$,
so \eqref{eq:kdv-quant-L2-0} follows from \eqref{eq:kdv-quant-L2}.
For \eqref{eq:kdv-quant-H2-0}, the proof of Theorem~\ref{thm:kdv-sobolev-quantitative}
applies verbatim with $q_n^g$ replaced by $p_{n,\varepsilon}^g$, with constants uniform in
$0<\varepsilon\le T$:
(i) Step~1 of Theorem~\ref{thm:order2} is proved first for the regularized solutions, so
$\sup_{n,t}\|p^g_{n,\varepsilon}(t)\|_{H^r}\le C_{r,T,R}\|g\|_{H^r}$ by \eqref{eq:high-reg-regbound};
(ii) in the two-mesh remainder the nonlinear mesh variation is
$|h_{n,\varepsilon}(t)^2-h_{m,\varepsilon}(t)^2|=\max\{t,\varepsilon\}^2|n^{-2}-m^{-2}|
\le T^2|n^{-2}-m^{-2}|$, while $E_n-E_m$ is unchanged;
(iii) Lemma~\ref{lem:H2onesided} is mesh-uniform, so the $H^2$ data-stability estimate
holds for $p_{n,\varepsilon}$;
(iv) $p_{m,\varepsilon}\to S(\cdot)g$ as $m\to\infty$ by \eqref{eq:fixed-epsilon-H3}.
Hence $\|p^g_{n,\varepsilon}-S(\cdot)g\|_{C_TH^2}\le C_{T,R}\,n^{-2}\|g\|_{H^7}$ for $g\in H^7$,
and Step~3 with $\delta=n^{-2/5}$ yields \eqref{eq:kdv-quant-H2-0}.}
\hfill$\square$
}

\subsection{$H^3$-energy estimate}
{
\begin{lemma}[Energy inequality for the  KdV flow]
\label{lem:kdv-global-H3-energy}
Let $u\in C([a,b];H^3(\R))\cap C^1([a,b];L^2(\R))$
solve $u_t=-D^3u+6uDu$. Then, for $a\le s\le t\le b$,
\begin{equation}\label{eq:kdv-H3-energy-ineq}
\|u(t)\|_{H^3}^2\le\|u(s)\|_{H^3}^2
+C\int_s^t\|Du(\tau)\|_\infty\|u(\tau)\|_{H^3}^2\,\dd\tau.
\end{equation}
\end{lemma}
\begin{proof}
Let $J_\eta$ be a standard Friedrichs convolution mollifier,
$u_\eta=J_\eta u$, and differentiate the $H^3$-energy 
$\sum_{k=0}^3\|D^ku_\eta\|_2^2$.
The Airy term has zero real energy pairing. For the nonlinear term,
write $D^k(uDu)=uD^{k+1}u+R_k(u)$. Then
Leibniz and Gagliardo--Nirenberg estimates yield
$\|R_k(u)\|_2\le C\|Du\|_\infty\|u\|_{H^3}$
for $k=0,\ldots,3$ (with $R_0=0$).
For $v=D^ku$, write
\[
J_\eta(uDv)=uD(J_\eta v)+[J_\eta,u]Dv.
\]
The first term is estimated by integration by parts and the uniform
Friedrichs commutator inequality
$\|[J_\eta,u]Dv\|_2\le C\|Du\|_\infty\|v\|_2$
controls the second. Hence
\[
\frac{\dd}{\dd t}\|u_\eta(t)\|_{H^3}^2
\le C\|Du(t)\|_\infty\|u(t)\|_{H^3}^2
\]
with $C$ independent of $\eta$; the analogous reversed inequality
is obtained by taking absolute values of the derivative.
Integrate on $[s,t]$ and send $\eta\downarrow0$, using
$J_\eta u(\tau)\to u(\tau)$ in $H^3$ at the fixed endpoints.
\end{proof}

}

\subsection{Proof of Theorem~\ref{prop:global-kdv}}
{
\emph{Proof of Theorem~\ref{prop:global-kdv}.}
We first work on an arbitrary compact subinterval of the maximal local
lifespan of $u=S(\cdot)f$. Consider
\begin{equation}\label{eq:global-invariants}
\begin{aligned}
\mathcal M(u)=\int_\R u^2\,\dd x,\qquad
\mathcal H(u)=\int_\R(\tfrac12u_x^2+u^3)\,\dd x,\qquad
\mathcal J(u)=\int_\R(u_{xx}^2+10uu_x^2+5u^4)\,\dd x.
\end{aligned}
\end{equation}
\emph{Conservation.}\footnote{See also \cite[Lemma~2.5]{Sjoberg1970}}
First let $f\in H^7$. 
The strong
$C_tH^2$ limit is then in $C_tH^6$ by Sobolev interpolation:
for $0\leq t\leq T$ and $n,m\geq 1$, 
\[
\|q_n(t)-q_m(t)\|_{H^6}
\le C\|q_n(t)-q_m(t)\|_{H^2}^{1/5}
       \|q_n(t)-q_m(t)\|_{H^7}^{4/5}\overset{\eqref{eq:two-mesh-cauchy},\eqref{eq:high-reg-propagation}}\leq C_{T,\|f\|_{H^7}}|n^{-2}-m^{-2}|^{1/5},
\]
The equation therefore gives $u\in C_t^1H^3$
Then,
\[
\frac{\dd}{\dd t}\int u_{xx}^2=30I,\quad
\frac{\dd}{\dd t}\int uu_x^2=-3I+6B,\quad
\frac{\dd}{\dd t}\int u^4=-12B,
\]
where $I=\int u_xu_{xx}^2$ and $B=\int uu_x^3$.
Hence, 
$\mathcal J(u(t))=\mathcal J(f)$. $\mathcal{M}$ and $\mathcal{H}$ are conserved by \eqref{eq:limit-mass} and \eqref{eq:KdV-fiels} respectively.

For general $f\in H^3$, approximate by
$f_j=\rho_{1/j}*f\in H^7$, using a Friedrichs mollifier.
Then $\|f_j\|_{H^3}\le\|f\|_{H^3}$ and $\|f_j\|_2\le\|f\|_2$,
so the common lifespan is not shortened. Corollary~\ref{cor:kdv-H2-continuous-dependence}
yields $S(\cdot)f_j\to S(\cdot)f$ in $C_tH^2$.
Since $\mathcal M,\mathcal H,\mathcal J:H^2\to\R$ are continuous,
the three invariants pass to the limit on that interval.
Uniqueness allows repetition on overlapping local intervals;
therefore all three invariants are conserved throughout the
maximal lifespan.

\emph{A priori bounds and continuation.}
The one-dimensional Gagliardo--Nirenberg inequality gives
$|\int u^3|\le C\|u\|_2^{5/2}\|Du\|_2^{1/2}$.
Conservation of $\mathcal M,\mathcal H$ and Young's inequality
therefore bound $\|u(t)\|_{H^1}$ uniformly throughout the lifespan.
Using $H^1(\R)\hookrightarrow L^\infty(\R)$ and conservation of
$\mathcal J$ then gives
\[
\|D^2u(t)\|_2^2
\le |\mathcal J(f)|+10\|u(t)\|_\infty\|Du(t)\|_2^2,
\]
where the nonnegative term $5\int u^4$ was discarded.
Consequently, $\|u(t)\|_{H^2}\le K_f$ on the entire maximal interval.
By Lemma~\ref{lem:kdv-global-H3-energy} and
$\|Du(t)\|_\infty\le C K_f$,
\begin{equation} \label{eq:H3-bnd0}
\|u(t)\|_{H^3}\le\|f\|_{H^3}
\exp(CK_f|t|).
\end{equation}

for each time in the lifespan (using reversibility for negative times).
If the forward maximal time $T_{\max}$ were finite, then
$\sup_{0\le t<T_{\max}}\|u(t)\|_{H^3}<\infty$.
The local lifespan depends only on the $H^3$ radius (and is bounded
from below on bounded subsets, by $\|f\|_2\le\|f\|_{H^3}$).
Restart the local construction from $u(t_0)$ with $t_0$ sufficiently
close to $T_{\max}$, yields a solution beyond $T_{\max}$.
Uniqueness in Lemma~\ref{lem:kdv-H2-uniqueness} identifies the two
solutions on their overlap, contradicting maximality.
Thus the solution extends forward for all time. By reversibility, the solution extends to $\R$. Strong $H^3$
continuity at each restart follows from
Proposition~\ref{prop:restart-H3}.

\emph{Stability and joint regularity on compact time intervals.}
Fix $A,\rho>0$. 
Then,
\[
M_{A,\rho}:=\sup_{\substack{|t|\le A\\\|f\|_{H^3}\le\rho}}
\|S(t)f\|_{H^3}\overset{\eqref{eq:H3-bnd0}}<\infty.
\]
Choose a positive local interval length $\tau=\tau(M_{A,\rho})$
from Corollary~\ref{cor:kdv-H2-continuous-dependence} and partition
$[0,A]$ into finitely many pieces of length at most $\tau$.
Restart the $H^2$-stability estimate on each piece; the
intermediate states remain in the same $H^3$ ball.
A finite product of the local Lipschitz constants proves
\eqref{eq:lipMain} for $0\le t\le A$,
 and reversibility yields it for
$-A\le t\le0$. Moreover, the KdV equation and the uniform bound
show $\|\partial_tS(t)f\|_2\le C(M_{A,\rho}+M_{A,\rho}^2)$.
Combining this time-Lipschitz bound in $L^2$ with
\eqref{eq:lipMain} proves \eqref{eq:global-joint-lipschitz}.
Finally, for each trajectory,
\[
\|S(t)g-S(s)g\|_{H^2}
\le C\|S(t)g-S(s)g\|_2^{1/3}
       \|S(t)g-S(s)g\|_{H^3}^{2/3}
\le C_{A,\rho}|t-s|^{1/3},
\]
Which together with \eqref{eq:lipMain} yields
\eqref{eq:global-time-holder}.
\hfill$\square$
}

\appendix
\section{Sobolev interpolation}\label{app:interpolation}
{
The following elementary interpolation estimate is useful when
only an $L^2$ error and a uniform higher-Sobolev bound are available.
}
\begin{lemma}[Sobolev interpolation]\label{lem:interp-H2}
For $0\le\sigma\le3$ and $f\in H^3(\R)$,
\begin{equation}\label{eq:interp-Hsigma}
 \norm[H^\sigma]{f}
 \le C_\sigma\norm[L^2]{f}^{1-\sigma/3}
       \norm[H^3]{f}^{\sigma/3}.
\end{equation}
In particular,
\begin{equation}\label{eq:interp-H2}
 {\displaystyle \norm[H^2]{f}\le C\norm[L^2]{f}^{1/3}
                         \norm[H^3]{f}^{2/3}.}
\end{equation}
\end{lemma}
\begin{proof}
It follows by Plancherel's theorem and H\"older's inequality.
\end{proof}

\section*{AI Disclosure.}\label{sec:ai-disclosure} The ideas underlying this preprint originated with the author. ChatGPT and Claude were used for mathematical discussions and editorial assistance. The author assumes full responsibility for the content.
{

\clearpage 
\bibliographystyle{unsrtnat}
\bibliography{kdv_references_exact}
}

\end{document}